\documentclass[10pt,a4paper]{article}
\makeatletter

\usepackage[left=3cm,top=3cm,bottom=3cm,right=3cm,head=1cm,foot=1cm]{geometry}
\usepackage{graphicx}
\usepackage{stmaryrd}
\usepackage{twoopt}
\usepackage{amssymb}
\usepackage{amsmath}
\usepackage{amsthm}
\usepackage{mathabx} 
\usepackage{amsfonts}
\usepackage{amscd}
\usepackage{url}
\usepackage[all]{xy}

\newcommand{\ch}{\mathrm{ch}}

\newcommand{\tr}{\mathrm{tr}}

\newcommand{\hil}{\mathcal{H}}

\newcommand{\g}{\mathfrak{g}}

\newcommand{\5}{\hspace{0,5cm}}
\newcommand{\3}{\vspace{0,3cm}}

\newcommand{\la}{\langle}
\newcommand{\ra}{\rangle}

\title{Twisted K-theory constructions in the case of a decomposable Dixmier-Douady class II: \\ Topological and Equivariant Models}
\date{}
\author{Antti J. Harju\footnote{University of Helsinki, email: harjuaj@gmail.com, tel.: +358 50 5724451}}
\begin{document}

\maketitle

\begin{abstract}
	This is a study of twisted K-theory on a product space $\mathbb{T} \times M$. The twisting comes from a decomposable cup product class which applies the 1-cohomology of $\mathbb{T}$ and the 2-cohomology of $M$. In the case of a topological product, we give a concrete realization for the gerbe associated to a cup product characteristic class and use this to realize twisted $K^1$-theory elements in terms of supercharge sections in a Fredholm bundle. The nontriviality of this construction is proved. Equivariant twisted K-theory and gerbes are studied in the product case as well. This part applies Lie groupoid theory. Superconnection formalism is used to provide a construction for characteristic polynomials which are used to extract information from the twisted K-theory classes. \3
	
\noindent MSC: 19L50, 53C08, 22A22.
\end{abstract}
\section*{Introduction}

In the first part of this series, \cite{HM12}, we gave a concrete realization for elements in a $K^1$-group twisted by a gerbe on a compact product manifold $\mathbb{T} \times M$. The twisting class decomposes as $\alpha \smile \beta \in H^3(\mathbb{T} \times M)$, for $\alpha \in H^1(\mathbb{T})$ and $\beta \in H^2(M)$. We used Hamiltonian quantization to construct the gerbes: quantization leads to a projective bundle of Hilbert spaces whose characteristic class can be read from the index of an odd Dirac family. The twisted K-theory classes were constructed explicitly using supersymmetric Fock spaces in the same spirit Wess-Zumino-Witten model has been used in the construction of equivariant twisted K-theory classes on a Lie group \cite{Mic04, FHT13}. One can also use projective families of Dirac operators to realize elements in twisted K-theory in the decomposable case, \cite{MMS05, MMS09, BG11}. 

A gerbe associated to the decomposable class has a simple topological realization. This is given by a bundle of projective Fock spaces over $\mathbb{T} \times M$.  The Fock space has a basis which can be decomposed into charge subspaces labeled by $\mathbb{Z}$. One fixes a line bundle $\lambda$ associated to a nontrivial class in $H^2(M, \mathbb{Z})$. Then one makes the Fock states of charge $k \in \mathbb{Z}$ to be topologically equivalent to $\lambda^{\otimes k}$ over the manifold $M$. Moreover, the translations over the circle $\mathbb{T}$ will raise the charge of each state. Consequently there are locally defined transition functions valued in the group of unitary operators on a Hilbert space. However, as one needs to create a state which is topologically a line bundle $\lambda$, these are naturally defined as local families of projective unitary transformations ($PU$-transformations). This bundle has a Dixmier-Douady class of the decomposable type: $\tau = \alpha \smile \beta$.  

Associated to these gerbes there are bundles of self-adjoint Fredholm operators acting on the Fock spaces. In \cite{HM12} we used supercharge operators in conformal field theory to construct sections in these bundles. The homotopy classes of supercharges determine elements in the twisted K-theory group $K^1(\mathbb{T} \times M, \alpha \smile \beta)$ which is isomorphic to an extension of 
\begin{eqnarray*}
	\{x \in K^1(M): x = x \otimes \lambda\} \5 \text{by} \5 \frac{K^{0}(M)}{K^{0}(M) \otimes (1 - \lambda)}.
\end{eqnarray*}
The supercharge sections are relevant to the subgroup isomorphic to $K^{0}(M) / K^{0}(M) \otimes (1 - \lambda)$. We used techniques from differential geometry to study the gerbes and supercharges in twisted $K$-theory. Especially superconnection techniques and families index theory were applied. Here we raise the issue of how to extend this formalism to the topological framework: the gerbes and twisted K-theory elements are still well defined on the topological product space $\mathbb{T} \times Y$. We give a complete characterization of the topological gerbes associated with the cup product Dixmier-Douady classes. The nontriviality of the supercharges as elements in twisted K-theory groups is a consequence of a spectral flow of the eigenvalues.  

We return to the smooth setup in the study of equivariant gerbes and equivariant twisted K-theory on a smooth manifold $\mathbb{T} \times M$. The group $G$ is a compact Lie group and acts smoothly on $M$. The action is extended to $\mathbb{T} \times M$ by letting $G$ act trivially on $\mathbb{T}$. The strategy to construct gerbes is to replace the twisting line bundles in the Fock bundle model with equivariant line bundles on $M$. We use the Lie groupoid formalism which is flexible enough to allow generalizations of this type. In this case we can give explicit formulas for the equivariant $PU$-transition functions for the gerbe in terms of cocycles in groupoid cohomology. The Dixmier-Douady class of the gerbe is the cup product class $\alpha \smile \beta_G$ where $\alpha$ is a 1-cohomology class associated to the trivial $G$-space structure in $\mathbb{T}$ whereas $\beta_G$ is a cocycle determined by a $G$-equivariant complex line bundle. 

Sections of $G$-invariant supercharge operators in a bundle of self-adjoint Fredholm operators are used to realize elements in the twisted equivariant K-theory groups on the product $\mathbb{T} \times M$. As in \cite{HM12} a Mayer-Vietoris sequence is used to solve the twisted K-theory group structure. We also develop equivariant characteristic polynomials for the homotopy classes of supercharges. This application uses equivariant superconnections. The target group of these characters is a quotient of the Cartan model of rational equivariant cohomology $H^{\mathrm{odd}}_G(\mathbb{T} \times M, \mathbb{Q})$. We need to go to a quotient group to make the characters periodic in $\mathbb{T}$. The target group for the superconnection is sensitive to some torsion coming from the twisted equivariant K-theory group. This target is essentially different from the twisted equivariant cohomology theory \cite{MS05}. 

As an example of equivariant gerbes and equivariant twisted K-theory we consider the action of the circle group $\mathbb{T}$ on the sphere $S^2$ which rotates the sphere around its axis. There are two fixed points: the poles. Consequently, there are nontrivial $\mathbb{T}$-equivariant line bundles parametrized by irreducible representations of $\mathbb{T}$ on the northern and the southern hemisphere. Such a line bundle is used to construct a $\mathbb{T}$-equivariant gerbe. The ordinary K-theory group $K^1_{\mathbb{T}}(\mathbb{T} \times M)$ is a sum of representation rings of $\mathbb{T}$ constrained by the common virtual dimension. The twisted $K^1$-groups are defined similarly but the generators of the representation rings are truncated by the choice of the twisting bundle $\lambda_{\mathbb{T}}$. Similar phenomenon occurs in twisted equivariant K-theory on Lie groups which was studied by Freed-Hopkins-Teleman, \cite{FHT11}, \cite{FHT13}. We apply the superconnection analysis to study these groups. 

\section{Twisted K-Theory of Topological Products}

\noindent \textbf{1.1.} Let $\hil$ denote an infinite dimensional complex Hilbert space. We use the model of Atiyah and Segal, \cite{AS04}, for the definition of twisted K-theory. One needs to modify the usual K-theoretic classification spaces since in the twisted case, the group of unitary transformations $U(\hil)$ needs to act continuously on them. Suppose that $U(\hil)$ is equipped with the strong operator topology and let $PU(\hil)$ denote the group of projective unitary transformations. $U(\hil)$ is now contractible. 

Let $\mathcal{B}(\hil)$ denote the space of bounded operators equipped with the strong operator topology. The space of compact operators $\mathcal{K}(\hil)$ is equipped with the norm topology. Let $\textbf{Fred}^{(0)}$ be the space of operators $A \in \mathcal{B}(\hil)$ which have a parametrix, an operator $B \in \mathcal{B}(\hil)$, such that $1 - AB$ and $1 - BA$ are compact. Let $\textbf{Fred}^{(1)}$ denote the subspace of self-adjoint operators with both positive and negative essential spectrum. 

Let $X$ be a compact topological space. Suppose that $\tau$ is a three cohomology class, $\tau \in H^3(X, \mathbb{Z})$ which, under the standard isomorphism $H^1(X, \underline{PU}(\hil)) \rightarrow H^3(X, \mathbb{Z}) $, is realized by a $PU(\hil)$-bundle $\textbf{P}$. Associated to the bundle $\textbf{P}$ we have a projective bundle of Hilbert spaces and the bundles of $\textbf{Fred}^{\bullet}$-operators for $\bullet = 0,1$. The group $PU(\hil)$ acts on the space of Fredholm operators under the adjoint action. The twisted K-theory groups are the groups of homotopy classes of continuous sections in the Fredholm bundles
\begin{eqnarray*}
K^{\bullet}(X, \tau) = \{[Q]: Q \in C(\textbf{P} \times_{\text{Ad}} \textbf{Fred}^{\bullet})\}.
\end{eqnarray*}
Precisely speaking, to define a twisted K-theory group one needs to fix a representative of the twisting class $H^3(X, \mathbb{Z})$. Namely, there is no canonical isomorphism relating two twisted K-theory groups with two different (but isomorphic) realizations of the twist class. However, in the analysis presented here, the projective bundle is always clearly defined and this should not lead to any confusion. \3

\noindent \textbf{1.2.} Let $\mathcal{H}$ be an infinite dimensional complex Hilbert space with a polarization $\mathcal{H} = \mathcal{H}^+ \oplus \mathcal{H}^-$. A canonical anti commutation relations algebra, CAR, is a unital $C^*$-algebra generated by $a(f): f \in \mathcal{H}$ such that $a: \mathcal{H} \rightarrow \text{CAR}$ is an antilinear map and CAR is subject to 
\begin{eqnarray*}
	\{a(u), a(v)\} = 0 \5 \text{and}\5 \{a(u), a(v)^* \} = \la u,v \ra \textbf{1},
\end{eqnarray*}
for all $u,v \in \mathcal{H}$. A Fock space $\mathcal{F}$ is a Hilbert space with a vacuum vector $|0 \rangle$ and CAR acts on the vacuum by $a(u)|0 \rangle = 0 = a^*(v)|0 \rangle$ for $u \in \mathcal{H}^+, v \in \mathcal{H}^-$ and such that the basis of $\mathcal{F}$ consists of 
\begin{eqnarray}\label{basis}
	a(u_{i_1}) \cdots a(u_{i_k}) a^*(u_{j_1}) \cdots a^*(u_{j_l}) |0 \rangle\5 \text{with} \5 u_{i_{x}} \in \mathcal{H}^-, u_{j_{y}} \in \mathcal{H}^+.  
\end{eqnarray}
This representation of CAR is irreducible. The basis vectors of this type are eigenvectors for the densely defined charge operator $N$, the eigenvalue of \eqref{basis} being $l-k$. Then we can write the Fock space as a completion of the charge subspaces $\mathcal{F} = \widehat{\bigoplus}_{k \in \mathbb{Z}} \mathcal{F}^{(k)}$. The operator $S$ denotes the shift operator which sends $\mathcal{F}^{(k)}$ to $\mathcal{F}^{(k+1)}$ for all $k \in \mathbb{Z}$. We refer to \cite{Ott95} for a careful introduction of these matters. 

The Hilbert space vectors will be labeled so that $\mathcal{H}^+$ has the basis $\{u_i : i \geq 0\}$ and $\mathcal{H}^-$ has the basis $\{u_i : i < 0\}$. \3

\noindent \textbf{1.3.}  Consider a compact topological space $Y$ and a nontrivial complex line bundle $\lambda$ on $Y$. We use projective Fock spaces in the construction of a gerbe over $\mathbb{T} \times Y$.

Let $\{V_i: 1 \leq i \leq n\}$ be an open cover of $Y$ which trivializes $\lambda$ with respect to the transition functions $h_{ij}: V_{ij} \rightarrow \mathbb{T}$. For $\mathbb{T}$ we choose a cover $\{\mathbb{T}_+, \mathbb{T}_-\}$ such that $\mathbb{T}_+ \cap \mathbb{T}_-$ consist of two disconnected arcs: neighborhoods of $-1$ and $1$ in $\mathbb{T}$.  We denote these intersections by $\mathbb{T}^{(\pm 1)}_{+ -}$. Concretely we can choose a small real $\epsilon > 0$ and then set $\mathbb{T}_+ = \{e^{i \phi}: \phi \in ( - \epsilon, \pi + \epsilon)\}$ and $\mathbb{T}_- = \{e^{i \phi}: \phi \in (\pi - \epsilon, 2 \pi + \epsilon)\}$. Whenever we consider the intersection$\mathbb{T}_{+-}^{(1)}$ we use $[\phi] = \phi \text{ mod } 2 \pi$. We define an open cover $\{U_i: 1 \leq i \leq 2n\}$ for $\mathbb{T} \times Y$ such that 
\begin{eqnarray*}
	U_i = \mathbb{T}_+ \times V_i \5 \text{and} \5 U_{i+n} = \mathbb{T}_- \times V_i\5 \hbox{for all} \5 1 \leq i \leq n. 
\end{eqnarray*}

For what follows we need a projective Fock bundle on $\mathbb{T} \times Y$ so that the translation around the circle $\mathbb{T}$ raises the charges of the Fock states in its fibres by one. This process will be localized in the open nbd $\mathbb{T}^{(+1)}_{+ - }$ of $\mathbb{T}$. More precisely, we write
\begin{eqnarray*}
U_{i,j+n}^{(1)} := \mathbb{T}^{(+1)}_{+ - } \times V_i \cap V_j 
\end{eqnarray*}
and  define the local families $S: U^{(1)}_{i,j+n} \rightarrow PU(\mathcal{H})$ for all $1 \leq i,j \leq n$ which act on the Fock states by raising the charge as in 1.2. These families are naturally valued in $PU(\mathcal{H})$: when acting on a Fock space, $S$ creates a state from the vacuum    
 \begin{eqnarray}\label{S}
	S |0 \rangle_{([\phi],p)} = a^*(u_0) |0 \rangle_{([\phi],p)}, \5 ([\phi] = \phi \text{ mod } 2 \pi)
\end{eqnarray}
but we need to choose an orthonormal basis vector $u_0$ from the subspace it spans. Thus, one needs to fix complex phases locally. We view $S:U^{(1)}_{i,j+n} \rightarrow PU(\mathcal{H})$ as components of transition functions of a $PU(\mathcal{H})$ bundle which relates fibres over $U_{j+n}$ and $U_i$ for $1 \leq i,j \leq n$. These operators act on the CAR algebra under conjugations by
 \begin{eqnarray*}
	S a(u_k) S^{-1} = a(u_{k+1}) 
\end{eqnarray*}
 for all $k \in \mathbb{Z}$ and similarly for $a^*(u_k)$. We can fix the basis vector $u_0$ and introduce lifted local families $\hat{S}: U_{i,j+n} \rightarrow U(\mathcal{H})$.  
 
Next we define a bundle of Fock spaces and apply the above local families to raise (lower) the charge of the fibres under translations in the positive (negative) direction in $\mathbb{T}$. In addition we make the charge $k$ states in the lifted Fock bundle to be of topological type $\lambda^{\otimes k}$ over $Y$. Such a projective bundle is constructed from the {\v C}ech-cocycle $g \in H^1(\mathbb{T} \times M, \underline{PU}(\mathcal{H}))$ whose nonidentity components are defined by 
\begin{eqnarray*}
	g_{ij}(\phi,p) &=& (h_{ij}(p))^N, \\ 
	g_{i+n,j+n}(\phi,p) &=& (h_{i+n,j+n}(p))^N, \\
	g^{(-1)}_{i,j+n}(\phi,p) &=& (h_{i,j+n}(p))^N \\
	g^{(-1)}_{j+n,i}(\phi,p) &=& (h_{j+n,i}(p))^N \\
	g^{(1)}_{i,j+n}([\phi],p) &=& (h_{i,j+n}(p))^{N} S  \\
	g^{(1)}_{j+n,i}([\phi],p) &=&  S^{-1} (h_{j+n,i}(p))^{N}
\end{eqnarray*}
The upper index $(\pm 1)$ means that the cocycle is defined in the neighborhood $\mathbb{T}^{(\pm 1)}_{+-}$ in $\mathbb{T}$. In particular, the domain of $g^{(1)}_{i,j+n}$ is the subset $U^{(1)}_{i,j+n}$.  Now we have the projective bundle $\textbf{PF}$ determined by the cocycle $\{g_{ij}: 1 \leq i,j \leq 2n\}$. Notice that the $\mathbb{T}$-valued components in this cocycle do not contribute at the level of projective transformations. However, when lifted, these become important. 

The gerbe defines a class in the {\v C}ech-cohomology group $H^2(\mathbb{T} \times Y, \underline{\mathbb{T}})$ whose components are determined by the lifted transition functions
\begin{eqnarray*}
	f_{ijk} = \hat{g}_{ij} \hat{g}_{jk} \hat{g}_{ik}^{-1}.
\end{eqnarray*}
for all $1 \leq i,j,k \leq 2n$. In \cite{HM12} we found out that the Dixmier-Douady cocycle is equal to the totally antisymmetric cocycle whose nonidentity components are determined by
\begin{eqnarray}\label{gerbe}
	f^{(1)}_{i,j+n,k+n}([\phi],p) = (h_{jk}(p))^{-1}
\end{eqnarray}
for $1 \leq i,j,k \leq n$. This is a consequence of the relation $SN = (N-1)S$ applied to the lifted transition functions. \3

\noindent \textbf{Proposition.} Under the isomorphism $H^2(\mathbb{T} \times Y, \underline{\mathbb{T}}) \rightarrow H^3(\mathbb{T} \times Y, \mathbb{Z})$, associated with the group extension $\mathbb{Z} \rightarrow \mathbb{R} \rightarrow \mathbb{T}$, the Dixmier-Douady class of $\textbf{PF}$ is represented by the cup product $\tau= \alpha \smile \beta$ where $\alpha$ is a generator of the group $H^1(\mathbb{T}, \mathbb{Z})$ and $\beta$ is the cohomology class of the twisting line bundle $\lambda$. \3

\noindent Proof. This can be proved by following the lines of the Proposition in 2.7 in the case of a trivial group $G = \{1\}$.\5 $\square$ \3

\noindent \textbf{1.4.} Here we recall the supercharge construction from \cite{HM12} which can be applied in a twisted K-theory of a topological space straightforwardly. Consider the real Clifford $*$-algebra generated by $\{ \psi_n, n \in \mathbb{Z} \}$ and subject to
\begin{eqnarray*}
	\{ \psi_n, \psi_m \} = 2 \delta_{n,-m}, \5 \psi_n^* = \psi_{-n}. 
\end{eqnarray*}
We use a vacuum representation, $\mathcal{H}_s$, for the Clifford algebra such that $\psi_0$ acts on the vacuum as the identity. Let $\eta_0$ denote the vacuum vector. The operators $\psi_i$ with $i < 0 $ annihilate the vacuum subspace whereas $\psi_i$ with $i>0$ are used to generate the basis of the module. 

Now we define a new bundle realization for the decomposable Dixmier-Douady class $\tau$. Let $\textbf{S}$ denote the canonically trivial bundle with the representations $\mathcal{H}_s$ as its fibres over $\mathbb{T} \times M$. The projective bundle $\textbf{PS}$ is a canonically trivial gerbe. Then we can define the $PU(\hil)$-bundle $\textbf{PS} \otimes \textbf{PF}$ using the tensor product structure of gerbes with a Dixmied-Douady class represented by the decomposable class $\tau$. The spinor bundle has been introduced so that we can produce Dirac type operators with both positive and negative essential spectrum. 

Associated to the gerbe $\textbf{PS} \otimes \textbf{PF}$ there is a bundle of bounded $\textbf{Fred}^{(1)}$-operators. Then, the homotopy classes of sections define the twisted K-theory groups $K^1(\mathbb{T} \times Y, \tau)$. We shall use the approximated sign map 
\begin{eqnarray*}
\Psi: A \mapsto \frac{A}{\sqrt{1+ A^2}} = F_A,
\end{eqnarray*}
to induce the topology from $\textbf{Fred}^{(1)}$ to the space of unbounded self-adjoint Fredholm operators $\textbf{Fred}^{(1)}_{\Psi}$. Associated to the gerbe $\textbf{PS} \otimes \textbf{PF}$ there is a bundle of $\textbf{Fred}^{(1)}_{\Psi}$-operators, and the homotopy classes of sections in this bundle define elements in $K^1(\mathbb{T} \times Y, \tau)$ through the correspondence $[Q] \mapsto [F_Q]$. 

The supercharges on the product manifolds are the local self-adjoint families $Q^i: U_i \rightarrow \textbf{Fred}_{\Psi}^{(1)}$ defined by 
\begin{eqnarray*}
	Q^i_{\phi,p} = \sum_{n \in \mathbb{Z}} \psi_n \otimes e_{-n} + \frac{\phi}{2\pi} \psi_0 \otimes \textbf{1}.
\end{eqnarray*}
We have used the operators $e_n$ which determine a representation of the centrally extended Lie algebra of the loop group $L \mathbb{T}$ on $\mathcal{F}$ by 
\begin{eqnarray*}
	e_n = \sum_{i \in \mathbb{Z}} : a^*(u_{n+i})a(u_i): \5 \text{with} \5 [e_n, e_m] = - n \delta_{n,-m}.
\end{eqnarray*}
The usual normal ordering :: is applied to make these operators well defined on the Fock space, $: a^*(u_{i})a(u_j): = - a(u_j)a^*(u_{i})$ if $ i = j < 0$ and the order is unchanged otherwise. Now $Q^i$ are families of unbounded and self-adjoint Fredholm operators. The $PU(\mathcal{H})$-valued cocycle $g$ acts on the local families $Q^i$ by conjugation according to  
\begin{eqnarray*}
    g_{i,j+n}([\phi],p) Q^{j+n}_{\phi,p} (g_{i,j+n}([\phi],p))^{-1} = Q^i_{\phi-2 \pi,p} 
\end{eqnarray*}
on each nonempty $U_{ij}$. Concretely, this is a consequence of the operator identities
\begin{eqnarray*}
S e_0 S^{-1} = e_0 - 1, \5 S e_n S^{-1} = e_n. 
\end{eqnarray*}
Therefore, $\{Q^i\}$ indeed defines a section in the unbounded Fredholm bundle. This construction is a straightforward generalization of \cite{HM12} to the topological setup. Let us write $F^i_Q(\phi,p) := F_{Q^i(\phi,p)}$. We conclude. \3

\noindent \textbf{Theorem.} The homotopy class of the section of bounded operators, $(\phi, p) \mapsto F^i_Q(\phi,p)$, defines an element in the twisted K-theory group $K^1(\mathbb{T} \times Y, \tau)$.\3

\noindent \textbf{1.5.}  Next we study conditions for nontriviality of the elements in twisted $K^1$-groups. Here we take $X$ to be any compact topological space with a nontrivial $3$-cohomology group. Let $\tau \in H^3(X, \mathbb{Z})$ and realize it as a projective $PU(\hil)$ principal bundle $\textbf{P}_{\tau}$. Consider an associated projective Hilbert bundle  $\textbf{PH}_{\tau}$ which restricts to a trivial bundle over a finite cover $\{U_i: i \in I\}$ of $X$. Denote by $g_{ij}$ the transition functions of the projective bundle and by $\hat{g}_{ij}$ the lifted transition functions. 

Below we shall work with local Hilbert bundles and lifted unitary transformations. Consider a section $Q$ of the associated bundle $\textbf{P}_{\tau} \times_{\text{Ad}} \textbf{Fred}^{(1)}_{\Psi}$. It can be viewed as local families of unbounded Fredholm operators $Q^i: U_i \rightarrow \textbf{Fred}^{(1)}_{\Psi}$ acting on the local Hilbert bundles over $U_i$. Let us denote by $\mathcal{H}_{i,x}^{(a, b)}$ the $(a,b)$-spectral subspaces of $Q^i_{x}$ in the fibre $\mathcal{H}_{i,x}$ of the local Hilbert bundle over $U_i$. More precisely, this subspace consists of the eigenvectors of $Q^i_x$ whose eigenvalues are in the open interval $(a,b)$.  On $U_{ij} \neq \emptyset$ these spectral subspaces transform as 
\begin{eqnarray*} 
\hat{g}_{ij} \mathcal{H}^{(a,b)}_{j,x} = \mathcal{H}^{(a,b)}_{i,x}.
\end{eqnarray*}
In this paper we are interested in the standard K-theoretic situation in which the twisted K-theory elements are approximated sign operators of unbounded Fredholm operators. The following assumption is standard: 
\begin{eqnarray}\label{Fredholm}
\text{for each $i \in I$, $x \in U_i$ and $a,b \in \mathbb{R}$, $\mathcal{H}_{i,x}^{(a, b)}$ are finite dimensional.}
\end{eqnarray}
A $\tau$-twisted spectral section $P$ of $Q$ is a section in the bundle $\textbf{P}_{\tau} \times_{\text{Ad}} \mathcal{B}(\hil)$ so that the values of $P$ are projection operators, and for some positive real number $R > 0$ and for every $i$ and $x \in U_i$ 
\begin{eqnarray}\label{spectralsection}
	P^i_{x} = \left\{ \begin{array}{ll} \textbf{1} & \text{ in } \mathcal{H}_{i,x}^{(R, \infty)},\\
	 0 & \text{ in } \mathcal{H}_{i,x}^{(- \infty,-R)}. \end{array} \right.
\end{eqnarray}

The proof of the following theorem applies the ideas of the corresponding result in \cite{MP97} in the ordinary K-theory. \3
 
\noindent \textbf{Theorem.} Let $Q$ be a section of a bundle of unbounded $\textbf{Fred}^{(1)}_{\Psi}$-operators so that the homotopy class of the approximated sign lies in a twisted K-theory group, $[F_Q] \in K^1(X, \tau)$. If \eqref{Fredholm} holds, then $Q$ has a $\tau$-twisted spectral section $\Leftrightarrow$ $F_Q$ represents the trivial element in $K^1(X, \tau)$. \3

\noindent Proof. Assume that the spectral section exists. Then define the local families $A^i_Q: U_i \rightarrow \mathcal{B}(\hil_i)$ 
\begin{eqnarray*}
	A_Q^i = \frac{Q^i}{\sqrt{(Q^i)^2 + 1}} - P^i \frac{Q^i}{\sqrt{(Q^i)^2 + 1}} P^i - (1 - P^i) \frac{Q^i}{\sqrt{(Q^i)^2 + 1}} (1-P^i).
\end{eqnarray*}
For each $i$ and $x \in U_i$, $A_Q^i(x)$ is the zero operator outside the spectral subspaces $\hil_{i,x}^{(-R,R)}$ which is finite dimensional. The local families $A_Q^i$ transform accordingly under $\text{Ad}(\hat{g}_{ij})$ because $P^i$ and $Q^i$ do.  Therefore, we have the equality of the twisted K-theory classes
\begin{eqnarray*}
[F_Q - A_Q] = [F_Q].
\end{eqnarray*}
Moreover, we have a homotopy connecting $F_Q - A_Q$ at $r = 0$ and  
\begin{eqnarray*}
	\widetilde{F}_Q^i(r) = P^i \frac{Q^i + r}{\sqrt{(Q^i)^2 + 1}}P^i + (1 - P^i)\frac{Q^i - r}{\sqrt{(Q^i)^2 + 1}}(1-P^i)
\end{eqnarray*}
for any $r \in \mathbb{R}$. If $r$ is large enough, the first term is strictly positive whereas the second term is strictly negative and thus $\widetilde{F}_Q(r)$ becomes invertible.

It remains to check that the homotopy $\widetilde{F}_Q^i(r)$ satisfies the continuity conditions in twisted K-theory for $|r| \leq K$ where $K$ is an arbitrarily large real number. The continuity of $\widetilde{F}_Q(r)$ in the strong operator topology of $\mathcal{B}(\hil)$ follows immediately from the continuity of $F_Q$, $P$ and $1 - F_Q^2 = (1 + Q^2)^{-1}$. According to \eqref{Fredholm}, the eigenvalues of $Q^i_x$ will go to infinity and all the eigenspaces must be finite dimensional. Therefore, the local families $r(1 + Q^2)^{-1/2}$ are valued in $\mathcal{K}(\hil)$ for $r$ arbitrary big. So, up to compact operators, $\widetilde{F}_Q^i(r)$ coincides with $P^i F_Q^i P^i + (1-P^i) F_Q^i (1-P^i)$ which coincides with $F^i_Q$ up to a finite rank operator. Therefore, we can take $\widetilde{F}_Q(r)$ to be its own parametrix.

It follows that $[F_Q] \in K^1(X, \tau)$ is a representative of the trivial element since $\widetilde{F}_Q(r)$ becomes invertible.\3

We use the following lemmas to prove the converse statement. \3

\noindent \textbf{Lemma} \cite{MP97}. Let $A$ be a bounded self-adjoint operator on $\mathcal{H}$ such that $A^2 - A$ is of finite rank. Assume that there is a self-adjoint projection $P$ such that $||A - P|| < 1/2$. Then $A: \mathrm{im}(P) \rightarrow \mathrm{im}(A P)$ is an isomorphism onto a closed subspace and if $\Pi$ denotes the orthogonal  projection onto $\mathrm{im}(A P)$, then $A - \Pi$ has a finite rank.\3 
	
\noindent \textbf{Lemma.} Let $B$ be a section of a $\tau$-twisted $\textbf{Fred}^{(1)}_{\Psi}$-bundle. If $U_{ij}\neq \emptyset$ and $P^i$ and $P^j$ are projections onto $\mathrm{im}(B^i)$ and $\mathrm{im}(B^j)$, then $\hat{g}_{ij} P^j (\hat{g}_{ij})^{-1} = P^i$ on $U_{ij}$. \3
	
\noindent Proof. Suppose that $x \in U_{ij} \neq \emptyset$. Then $ \hat{g}_{ij}(x) B^j_x \hat{g}_{ij}(x)^{-1} = B^i_x$ and so 
\begin{eqnarray*}
 \text{im}(\hat{g}_{ij}(x) B^j_x \hat{g}_{ij}(x)^{-1}) = \text{im}(B^i_x). 
\end{eqnarray*}
Thus, the projections have to satisfy $\hat{g}_{ij}(x) P^j_x \hat{g}_{ij}(x)^{-1} = P^i_x$. \3

Assume that $F_Q$ represents the trivial twisted K-theory element. Then, at the level of unbounded operators, we can set a homotopy $Q(t)$ in the space of sections  starting with a section of invertible self-adjoint operators $Q(0)$ and ending with $Q(1) = Q$. In addition, we can choose $Q(0)$ such that there are no eigenvalues in $(0, \epsilon)$ where $\epsilon$ is a small positive real number. 

We define a continuous function $\chi \in C(\mathbb{R})$ such that $\chi(\lambda) = 0$ for $\lambda <0$ and $\chi(\lambda) = 1 $ for $\lambda \geq \epsilon$ and apply the spectral theorem to define
\begin{eqnarray*}
	J^i(t) = \chi(Q^i(t)), \5 t \in I:= [0,1].
\end{eqnarray*}
Then $J(0)$ is a section of projection operators and $J(1)$ satisfies the condition \eqref{spectralsection} for $R = \epsilon$ but $J(1)$ may not be a projection section, however, $(J(1))^2 -J(1)$ is a finite rank section.

The strategy of the proof is to divide the interval $[0,1]$  to subintervals $[t',t'']$ on which 
\begin{eqnarray*}
|| J^i_{x}(t') - J^i_{x}(s) || < \frac{1}{2}, \5 \text{for all} \5 s \in [t',t'']. 
\end{eqnarray*}
Then we proceed inductively by proving that if the lower limit in the interval has a projection section, $P(t')$, such that $J(t') - P(t')$ is a finite rank section, then such a section exists for any $t$-value in this interval. 

At $t = 0$ a projection section exists. Consider $s \in [t',t'']$ and define another section of $\mathcal{B}(\hil)$-operators 
\begin{eqnarray*}
A^i_{x}(s) = J^i_{x}(s) + P^i_{x}(t') - J^i_{x}(t').
\end{eqnarray*}
Now $||A^i_{x}(s) - P^i_{x}(t')|| < 1/2$ and $A^i_{x}(s)$ is equal to $P^i_{x}(t')$ up to a finite rank operator which implies that $A^i_{x}(s)^2 - A^i_{x}(s)$ is a finite rank operator for all $i$ and $x \in U_i$. By the first lemma above, the orthogonal projection $P^i_{x}(s)$ onto the subspace $\mathrm{im}(A^i_{x}(s) P^i_{x}(t'))$ in the fibre at $x$ of the local Hilbert bundle over $U_i$ differs from $A^i_{x}(s)$ and thus from $J^i_{x}(s)$ by a finite rank operator. By the second lemma, we can patch together these local projection families to form a section of $\mathcal{B}(\hil)$-valued operators. Thus, $ t \mapsto P(t)$ defines a homotopy of projection sections and $J(t) - P(t)$ is a finite rank section for all $t \in [0,1]$.

Consider the decomposition of the fibres $\mathcal{H}_{i,x}$ in the local Hilbert bundles into the eigenspaces of $Q^i$. Then $J(1)$ is diagonal on the fibres. Next we construct a new section $\widetilde{P}^i$ of projective operators so that $J^i_x(1) - \widetilde{P}^i_x$ will be an operator which is nonzero only in a subspace spanned by a finite number of $Q^i_x$ eigenspaces for all $x$ and $i$. Given such $\widetilde{P}$, we can take $R$ to be a positive real so that this subspace will lie in $\mathcal{H}^{[-R,R]}_{i,x}$ for all $i$ and $x$. We can also make $R$ to satisfy $R \geq \epsilon$. Then $\widetilde{P}$ is a twisted spectral section.

Let $J(1) - P(1) = C$. Then we continuously truncate the $Q$-eigenspace expansion of $C$ to produce a new section $\hat{C}$ such that the image of $x \mapsto \hat{C}^i_x$ is a subspace spanned by a finite number of eigenvectors of $Q^i_x$. Let $\widetilde{P}$ be the projection onto the image of
\begin{eqnarray*}
(J(1) - \hat{C})P(1) = (J(1) - \hat{C})(J(1) - C).
\end{eqnarray*}
By studying the $Q$-eigenspace decompositions of this sections, one easily verifies that $\widetilde{P}^i_x$ have the required properties for all $i$ and $x \in X$. Thus, we have found a twisted spectral section $\widetilde{P}$. \5 $\square$ \3

As was observed in \cite{MP97}, the spectral flow around the circle is an obstruction to define a spectral section. This is because all the eigenvalues will grow by $1$ leading to identifications of eigenstates corresponding to the eigenvalues $k$ and $k + 1$ at $\phi = 2 \pi$. The twisted spectral section would cut the eigenspace decomposition at some point which makes it noncontinuous. \3

\noindent \textbf{Corollary.} The supercharge sections $Q$ on $\mathbb{T} \times Y$ define a nontrivial element in $K^1(\mathbb{T} \times Y, \tau)$. \3

In \cite{HM12} we studied these families on compact manifolds and introduced another dependency on the topology of $M$ by tensoring the gerbe with an arbitrary rank complex line bundle $\xi$. We proved an index theorem using the superconnection techniques and observed that the nontrivial twisted K-theory information is localized around the null space of $Q$ and that the twisted K-theory class depends on the Chern character of $\xi$. Here the same phenomenon is expected, although the families index formulas are not available in the purely topological setup.

\section{Twisted K-theory of Action Groupoids}

\noindent \textbf{2.1.}  Equivariant twisted K-theory is convenient to describe using the theory of groupoids. Here we return to the smooth setup and study equivariant product manifolds in the language of Lie groupoids. Suppose that $\Gamma \rightrightarrows X$ is a Lie groupoid and $X$ is a compact manifold. 

For each $p \geq 1$ let $\Gamma^{(p)}$ denote the manifold of composable sequences of arrows. We also use $\Gamma^{(0)} = X$. We have $p+1$ maps 
 \begin{eqnarray*}
	\partial_i: \Gamma^{(p)} \rightarrow \Gamma^{(p-1)}\5 \text{for}\5 0 \leq i \leq p
\end{eqnarray*}
 such that $\partial_0$ leaves out the first arrow, $\partial_p$ leaves out the last arrow and the rest of $\partial_i$ compose the $i$'th and $i+1$'th arrow. Then we set $\partial: \Gamma^{(p)} \rightarrow  \Gamma^{(p-1)}$ by $\partial = \sum_{i} (-1)^{i}\partial_i$.
 
Let $F$ denote an abelian sheaf in the category of smooth manifolds and $F^{p}$ a small sheaf induced by $F$ on $\Gamma^{(p)}$.  Then we choose injective resolutions $F^{p} \rightarrow I^{p \bullet}$.  This defines a double complex $I^{\bullet}(\Gamma^{\bullet})$ whose total cohomology groups, $H^k(\Gamma, F)$, are the cohomology groups of the groupoids with values in the abelian sheaf $F$. This theory is Morita invariant \cite{BX11}, \cite{TXL04}. 

We can also use the de Rham complex in the case of coefficients in $\mathbb{R}$. The exterior derivative defines a map $d: \Lambda^k(\Gamma^{(p)}) \rightarrow \Lambda^{k+1}(\Gamma^{(p)})$. We can use this with the coboundary operator $\partial^*: \Lambda^k(\Gamma^{(p)}) \rightarrow  \Lambda^{k}(\Gamma^{(p+1)})$ to construct a double complex. The total cohomology of this complex is the de Rham cohomology of $\Gamma$, we write $H_{\mathrm{dR}}^k(\Gamma) = H^k(\Lambda^*(\Gamma^{\bullet}))$. The de Rham theorem is still valid in the form 
\begin{eqnarray*}
	H^k_{\mathrm{dR}}(\Gamma) \simeq H^k(\Gamma, \mathbb{R}).
\end{eqnarray*}
A cocycle in the de Rham theory is an integer cocycle if it lies in the image of the canonical map $H^k(\Gamma, \mathbb{Z}) \rightarrow H^k_{\text{dR}}(\Gamma, \mathbb{R})$. 
 
 We are interested in the cohomology groups $H^2(\Gamma, \underline{\mathbb{T}})$ and $H^3(\Gamma, \underline{\mathbb{Z}})$. We also point out the existence of the group $H^1(\Gamma, \underline{PU}(\mathcal{H}))$. We only consider compact groupoids (and therefore proper) and consequently the group homomorphism $H^2(\Gamma, \underline{\mathbb{T}}) \rightarrow H^3(\Gamma, \underline{\mathbb{Z}})$, constructed from the standard exact sequence of groups 
\begin{eqnarray*}
	0 \rightarrow \mathbb{Z} \stackrel{2 \pi i}{\longrightarrow} i\mathbb{R} \stackrel{\exp}{\longrightarrow} \mathbb{T} \rightarrow 1,
\end{eqnarray*} 
  is an isomorphism. The group homomorphism $H^1(\Gamma, \underline{PU}(\mathcal{H})) \rightarrow H^2(\Gamma, \underline{\mathbb{T}})$, constructed from the central extension $\mathbb{T} \rightarrow U(\mathcal{H}) \rightarrow PU(\mathcal{H})$, is not an isomorphism in general. However, there is a canonical left inverse, \cite{TXL04}. 
 
 A $\mathbb{T}$-central extension of a groupoid $\Gamma \rightrightarrows X$ is a groupoid $\hat{\Gamma} \rightrightarrows X$ together with a
groupoid morphism
\begin{eqnarray*}
(\pi, \mathrm{id}): (\hat{\Gamma} \rightrightarrows X) \rightarrow (\Gamma \rightrightarrows X)
\end{eqnarray*}
such that there is a left $\mathbb{T}$ action on $\hat{\Gamma}$ making $\pi: \hat{\Gamma} \rightarrow \Gamma$ a (left) principal
$\mathbb{T}$-bundle. These two structures are required to satisfy the
compatibility
\begin{eqnarray*}
(\mu \cdot x)(\mu' \cdot y) = \mu \mu' \cdot (xy)
\end{eqnarray*}
for all $(x,y) \in \hat{\Gamma}^{(2)}$ and $\mu, \mu' \in \mathbb{T}$. 

In \cite{TXL04} the elements of $\text{Ext}(\Gamma,\mathbb{T})$, the Morita equivalence classes of $\mathbb{T}$ extensions of $\Gamma$, were identified with isomorphism classes of $\mathbb{T}$-gerbes over the stack associated to $\Gamma$. The set $\text{Ext}(\Gamma,\mathbb{T})$ can be equipped with a structure of an abelian group. Then the group $\text{Ext}(\Gamma,\mathbb{T})$ is isomorphic to $H^2(\Gamma, \underline{\mathbb{T}})$. We shall not need the notion of an gerbe over a stack or Morita equivalence in this work and we refer to \cite{BX11}, \cite{TXL04} for details. In our case, the twisting classes arise naturally from equivariant $PU(\mathcal{H})$-bundles. Therefore, we are given an $H^1(\Gamma, \underline{PU}(\mathcal{H}))$ cocycle which then determines a class in $H^2(\Gamma,\underline{\mathbb{T}})$. This is related to the decomposable class $\alpha \smile \beta_G$ under the isomorphism onto $H^3(\Gamma, \mathbb{Z})$. \3

\noindent \textbf{2.2.} In this application we start with a cocycle $g \in H^1(\Gamma, \underline{PU}(\mathcal{H}))$. The lifting problem for the projection $U(\hil) \rightarrow PU(\hil)$ defines a class in $H^2(\Gamma,  \underline{\mathbb{T}})$. The cocycle $g$ defines a principal $PU(\mathcal{H})$ bundle $\textbf{P} \rightarrow X$ over the groupoid $\Gamma \rightrightarrows X$. Associated to $\textbf{P}$ there is a bundle of Fredholm operators
\begin{eqnarray*}
\textbf{P} \times_{\text{Ad}} \textbf{Fred}^{(1)} \rightarrow X.
\end{eqnarray*}
We denote by $C_{\Gamma}(\textbf{P} \times_{\text{Ad}} \textbf{Fred}^{(1)})$ the space of continuous $\Gamma$-invariant sections of the bundle $\textbf{P} \times_{\text{Ad}} \textbf{Fred}^{(1)}$. Recall from 1.1 that continuity of a section $A$ getting values in $ \textbf{Fred}^{(1)}$ means that there is a parametrix section $B$ such that $1 - AB$ and $1 - BA$ are continuous compact operators and both $A$ and $B$ are strongly continuous sections of Fredholm operators. The twisted $K^1$-group on $\Gamma$ is defined by  
\begin{eqnarray*}
	K^1(\Gamma, \tau) = \{ [Q] : Q \in C_{\Gamma}(\textbf{P} \times_{\text{Ad}} \textbf{Fred}^{(1)}) \}
\end{eqnarray*}
if $\tau$ is the twisting cohomology class and $\textbf{P}$ is its bundle theoretic realization.

As above, we will work with unbounded operators and the twisted K-theory elements will be realized by mapping them to the space of bounded operators via the approximated sign construction.\3

\noindent \textbf{2.3.} Suppose that $G$ is a compact Lie group acting smoothly on $M$. Here we consider finite groups as compact Lie groups of dimension zero. The goal is to give an explicit description for the equivariant gerbe in terms of transition functions over a Lie groupoid representing the equivariant manifold. We shall approach as follows. We start with the action groupoid $G \ltimes M \rightrightarrows M$ and associated to it define a {\v C}ech-type groupoid with local group actions. Such a groupoid is Morita equivalent to the action groupoid provided that we can pick a G-invariant cover. In addition, we shall apply the groupoid double complex cohomology associated with a {\v C}ech resolution, and for this reason the following condition is introduced:
\begin{quote}
There exists a $G$-invariant cover $\mathfrak{V} = \{V_{i}: 1 \leq i \leq n \}$ of $M$ such that all intersections of its components, $V_{i_1} \cap \cdots \cap V_{i_l}$ with $l \geq 1$, have trivial $\underline{\mathbb{T}}$-valued {\v C}ech cohomology in nonzero degrees:    
\begin{eqnarray}\label{Cech}
H^k(V_{i_1} \cap \cdots \cap V_{i_l} ,\underline{\mathbb{T}}) = 0 \5 \text{for all $k \geq 1$}.
\end{eqnarray}
\end{quote}

This condition holds if we have a $G$-equivariant good cover. However, it is essentially weaker. In 3.3 we consider a case with a $G$-equivariant cover satisfying the condition \eqref{Cech} which is not good.\3

\noindent \textbf{2.4.} Fix a $G$-invariant cover $\{V_i: 1 \leq i \leq n\}$ of $M$. We apply the trivial $G$-action on the unit circle. Then the cover $\{\mathbb{T}_{\pm} \times V_i: 1 \leq i \leq n\}$ is $G$-invariant. We shall denote this cover by $\{U_i : 1 \leq i \leq 2n\}$ as in 1.1. Consider the Lie groupoid 
\begin{eqnarray*}
\Gamma := \coprod_{ij} G \times  U_{ij} \rightrightarrows \coprod_i U_{i}
\end{eqnarray*}
with $s(g,x_{ij}) = x_j$ and $t(g, x_{ij}) = g x_i$. The multiplication of composable pairs is defined by 
\begin{eqnarray*}
	(g,(hx)_{ij}) (h, x_{jk}) = (gh, x_{ik}).
\end{eqnarray*}
From \cite{TXL04} and \cite{BX11} we get the following technical result. \3

\noindent \textbf{Proposition.}
\begin{quote}
\noindent (1) The groupoid $\Gamma$ is Morita equivalent to the action groupoid $G \ltimes (\mathbb{T} \times M)$.

\noindent (2) The Lie groupoid cohomology groups $H^*(\Gamma, \mathbb{Z})$ are isomorphic to the equivariant cohomology groups $H^*_G(\mathbb{T} \times M)$ with integer coefficients.
\end{quote}

Since $M$ is compact, the equivariant cohomology can be computed with the Cartan or Weyl model. We shall use the Cartan model in 2.11 to map the twisted K-theory elements to equivariant cohomology.\3

\noindent \textbf{2.5.} Here we forget $\mathbb{T}$ for a moment and consider $G$-equivariant complex line bundles over $M$. They are classified, up to an isomorphism, by the equivariant cohomology group $H^2_G(M)$. Therefore they are also classified by the 2-cohomology of the Morita equivalent action Lie groupoid
\begin{eqnarray*}
 \Xi := \coprod_{ij} G \times V_{ij} \rightrightarrows \coprod_i V_i,
\end{eqnarray*}
or by the isomorphic sheaf cohomology group $H^1(\Xi, \underline{\mathbb{T}})$. 

Recall that in the case of a compact manifold, the sheaf cohomology associated to an injective resolution is isomorphic to a {\v C}ech cohomology computed from an open cover which is acyclic for the sheaf. We shall consider the sheaf $\underline{\mathbb{T}}$. According to \eqref{Cech}, $\Xi^{(0)} = \coprod V_i$ can be covered with just one set, $\Xi^{(0)}$ itself. For the space of arrows $\Xi^{(1)}$ we need to fix a good cover $\{N_a: a \in J\}$ of $G$ which makes the complex acyclic. Then we fix a cover for $\Xi^{(1)}$ by 
\begin{eqnarray*}
\{\Xi^{(1)}_a = \coprod_{ij} N_a \times V_{ij}, a \in J\}.
\end{eqnarray*}
This cover makes the {\v C}ech resolution acyclic. For $\Xi^{(i)}$ with $i \geq 2$, one proceeds similarly and fixes a good cover in the direction of the group elements and applies the {\v C}ech construction. 

The groupoid double complex is equivalent to the {\v C}ech hypercohomology complex 
\begin{eqnarray}\label{double}
	\xymatrix{&\cdots & \cdots & \cdots \\ & C^{\infty}(\Xi^{(0)}, \underline{\mathbb{T}}) \ar[u]^{\delta} \ar[r] & \prod_{ab}C^{\infty}(\Xi_{ab}^{(1)}, \underline{\mathbb{T}}) \ar[u]^{\delta} \ar[r] &  \prod_{ab}C^{\infty}(\Xi_{ab}^{(2)},\underline{\mathbb{T}}) \ar[u]^{\delta} \ar[r] & \cdots \\ & C^{\infty}(\Xi^{(0)},\underline{\mathbb{T}})\ar[u]^{\delta} \ar[r] & \prod_{a}C^{\infty}(\Xi_{a}^{(1)},\underline{\mathbb{T}})\ar[u]^{\delta} \ar[r] &  \prod_{a}C^{\infty}(\Xi_{a}^{(2)},\underline{\mathbb{T}}) \ar[u]^{\delta} \ar[r] & \cdots } 
\end{eqnarray}
The horizontal lines are induced from the isomorphisms to the injective resolution and the choice is not unique. However, the cohomology groups associated to different choices are isomorphic. 

A 1-cocycle in this double complex has components $h^{01}$ and $h^{10}$ (the first index labels the vertical direction). However, the second vertical arrow in the first column is the identity map and consequently, the only nontrivial component is $h^{01} \in \prod_a C^{\infty}(\Xi_a^{(1)},\underline{\mathbb{T}})$. The cocycle condition in the vertical direction requires that $\delta(h^{01}) = 1$ and therefore $h^{01}$ can be glued over the components $N_a$ to a a sequence of functions $h_{ij}: \Xi^{(1)} \rightarrow \mathbb{T}$ which we denote by $h$. Also, $h$ is in the image of the injective map from the groupoid cohomology complex $(\Xi^{\bullet}, \partial^*)$ to the double complex and the cocycle condition reads $\partial^*(h) = 0$ which is equivalent to 
\begin{eqnarray}\label{linebundle}
	h_{ij}(g,hx) h_{jk}(h,x) = h_{ik}(gh,x)
\end{eqnarray}
on $V_{ijk}$. This analysis shows that it is sufficient to study the groupoid cohomology complex $(\Xi^{\bullet}, \partial^*)$ to realize the $G$-equivariant line bundles: any $G$-equivariant complex line bundle on $M$ is determined by a set of transition functions $\{h_{ij}: V_{ij} \rightarrow \mathbb{T}\}$ subject to $\partial^*(h) = 0$. The $G$-action on the fibres is given by 
\begin{eqnarray*}
(x_j, v) \mapsto ((ax)_i, h_{ij}(a,x)v)
\end{eqnarray*}
where $a \in G$, $(x_j, v)$ is written in the trivialization over $U_j$ and  $((ax)_i, h_{ij}(a,x)v)$ over $U_i$.

We can use the standard projection $p: \Gamma^{(1)} \rightarrow \Xi^{(1)}$ to pull the $G$-equivariant line bundle data from $M$ to $\mathbb{T} \times M$. \3 

\noindent \textbf{2.6.} Next we define a $G$-equivariant projective Fock bundle over $\mathbb{T} \times M$ in the spirit of \cite{HM12} or 1.1. The twisting is now determined by a $G$-equivariant line bundle $\lambda_G$ on $M$.  We denote by $h \in H^1(\Gamma, \underline{\mathbb{T}})$ a representative of this bundle in the Lie groupoid cohomology. 

Consider the operators $S$ and $N$ acting on the Fock space as in 1.2. Again, $S$ is defined up to fixing of a basis for the state it creates and therefore it is naturally defined as a $PU(\mathcal{H})$ operator. Then we define a projective Fock bundle, $\textbf{PF}_G$, associated to the cocycle  $g \in H^1(\Gamma, \underline{PU}(\mathcal{H}))$ with the nonidentity components $g_{ij}: G \times U_{ij} \rightarrow PU(\hil)$ given by
 \begin{eqnarray*}
	g_{ij}(a, x) &=& (h_{ij}(a,x))^N \\ 
	g_{i+n,j+n}(a, x) &=& (h_{i+n,j+n}(a,x))^N \\
	g^{(-1)}_{i,j+n}(a,x) &=& (h_{i,j+n}(a,x))^N \\
	g^{(-1)}_{j+n,i}(a,x) &=& (h_{j+n,i}(a,x))^N \\
	g^{(1)}_{i,j+n}(a,x) &=& (h_{i,j+n}(a,x))^{N} S \\
	g^{(1)}_{j+n,i}(a,x) &=&  S^{-1}(h_{j+n,i}(a,x))^{N}
\end{eqnarray*}
 for all $1 \leq i,j \leq n$.  The translation around the circle in the positive direction adds the charge by one and the translation in the negative direction subtracts the charge by one. In local coordinates the action of $G$ reads
\begin{eqnarray*}
	(x_j, \Psi) \mapsto ((ax)_i, g_{ij}(a,x) \Psi) \5 \hbox{for $1 \leq i,j \leq 2n$}
\end{eqnarray*}
if $(x_j, \Psi)$ is evaluated in the trivialization over $U_j$ and $((ax)_i, g_{ij}(a,x) \Psi)$ in the trivialization over $U_i$ and $\Psi$ is a vector in the Fock space.

Denote by $\hat{g} \in H^1(\Gamma, \underline{U}(\hil))$ the lifted transition functions. They act on the Fock states over $M$ so that each charge $k$ subbundle has a structure of a $G$-equivariant vector bundle of type $\lambda_G^{\otimes k}$. Over $\mathbb{T} \times M$ there is a $G$-equivariant gerbe. The components of the Dixmier-Douady class $f \in H^2(\Gamma, \underline{\mathbb{T}})$ are determined by the formulas
\begin{eqnarray*}
    f_{ijk}&:& \Gamma^{(2)} \rightarrow \mathbb{T}, \\
    f_{ijk}(a,bx,b,x) &=& \hat{g}_{ij}(a,bx) \hat{g}_{jk}(b,x) (\hat{g}_{ik}(ab,x))^{-1}
\end{eqnarray*}
for the lifted transition functions. The following is a direct computation. \3

\noindent \textbf{Proposition.} The nonidentity components of the cocycle $f \in  H^2(\Gamma, \underline{\mathbb{T}})$ are
\begin{eqnarray*}
   f^{(1)}_{i,j+n,k+n}(a,bx,b,x) &=& (h_{j+n,k+n}(b,x))^{-1}  \\
   f^{(1)}_{j+n,i,k+n}(a,bx,b,x) &=& (h_{j+n,k+n}(ab,x)) \\
   f^{(1)}_{j+n,k+n,i}(a,bx,b,x) &=& (h_{j+n,k+n}(a,bx))^{-1}
\end{eqnarray*}
for $1 \leq i,j,k \leq n$ and $a,b \in G$. \3

\noindent \textbf{2.7.} Let $\lambda_G$ and $h \in H^1(\Gamma, \underline{\mathbb{T}})$ denote a $G$-equivariant complex line bundle and its structure cocycle as in 2.5. Let $\beta_G$ denote the cocycle in $H^2(\Gamma, \mathbb{Z})$ which is the image of $h$ under the isomorphism of cohomology groups (recall 2.1). Define $\alpha \in H^1(\Gamma,\mathbb{Z})$ to be the groupoid cocycle with the nonzero components
\begin{eqnarray*}
	\alpha_{i,j+n}^{(1)}(a,x) = 1 \5 \text{and} \5 \alpha_{i+n,j}^{(1)}(a,x) = -1
\end{eqnarray*}
for all $1 \leq i,j \leq n$ and $(a,x) \in G \times (\mathbb{T} \times M)$. We can view $\alpha$ and $\beta_G$ as groupoid de Rham cocycles. Then we can apply the cup product in $H^*_{\mathrm{dR}}(\Gamma, \mathbb{Z})$, \cite{BX11}. \3

\noindent \textbf{Proposition.} The Dixmier-Douady class $f \in H^2(\Gamma, \mathbb{Z})$ of the equivariant gerbe maps to the cup product class $\alpha \smile \beta_G$ under the group isomorphism $H^2(\Gamma, \underline{\mathbb{T}}) \rightarrow H^3(\Gamma, \mathbb{Z})$ induced by the standard group extension sequence 
\begin{eqnarray*}
\mathbb{Z} \stackrel{2 \pi i}{\longrightarrow} i \mathbb{R} \stackrel{\exp}{\longrightarrow} \mathbb{T}. 
\end{eqnarray*}

\noindent Proof. Denote by $h' \in H^1(\Gamma, \underline{i \mathbb{R}})$ the cocycle such that $(2 \pi i)^{-1} \partial^*(h') = \beta_G$ and $\exp(h') = h$. The nonidentity components of $\alpha \smile \beta_G$ are 
\begin{eqnarray*}
	&&(\alpha \smile \beta)^{(1)}_{i,j+n,k,l}(a, b c x, b, c x, c, x) = \frac{1}{2\pi i}(\partial^* h')^{(1)}_{j+n,k,l}( b, c x, c, x) \\
&&\5 \text{and} \\
   &&(\alpha \smile \beta)^{(1)}_{i+n,j,k,l}(a, b c x, b, c x, c, x) = - \frac{1}{2\pi i}(\partial^* h')^{(1)}_{jkl}( b, c x, c, x)
\end{eqnarray*}
for $1 \leq i,j \leq n$ and $1 \leq k,l \leq 2n$. The cup product $\alpha \smile \beta_G$ is equal to $\frac{1}{2\pi i} (\partial^*t)$ where $t \in H^2(\Gamma^{\bullet}, i \underline{\mathbb{R}})$ is the cocycle with the nonzero components 
\begin{eqnarray*}
	&&t^{(1)}_{i,j+n,k}(a, b x, b,x) = - (h')^{(1)}_{jk}(b, x) \5 \text{and} \\ 
	&&t^{(1)}_{i+n,j,k}(a, b x, b,x) = (h')^{(1)}_{jk}(b, x) 
\end{eqnarray*}
for $1 \leq i,j \leq n$ and $1 \leq k \leq 2n$. The image of the exponential cocycle $\exp(t)$ is the the cup product cocycle under the isomorphism $H^2(\Gamma, \underline{\mathbb{T}}) \rightarrow H^3(\Gamma, \mathbb{Z}) $. 

We introduce another cochain, $s \in C^1(\Gamma, \underline{\mathbb{T}})$ whose nonzero components are 
\begin{eqnarray*}
	s^{(1)}_{i+n,j}(a,x) = h_{i+n,j}(a,x) 
\end{eqnarray*}
for all $1 \leq i,j \leq n$. Then one checks that 
\begin{eqnarray*}
	f = \partial^*(s) \exp(t)
\end{eqnarray*}
It follows that the cohomology class of $f$ maps to the cohomology class of $\alpha \smile \beta_G$ under the isomorphism. \5 $\square$ \3

The gerbe $\textbf{PF}_G$ is nontrivial since its Dixmier-Douady class can be realized as a cup product of two nontrivial cocycles which are localized on different manifolds. Moreover, we note that if we cut the circle $\mathbb{T}$, then the Dixmier-Douady class becomes trivial since the 1-cocycle trivializes. \3

\noindent \textbf{2.8.} The twisted K-theory groups $K^*(\Gamma, \tau)$ associated to a gerbe with a Dixmier-Douady class represented by $\tau = \alpha \smile \beta_G$ can be solved, up to a group extension problem, using the Mayer-Vietoris sequence. Let us fix a closed cover $\mathbb{T} \times M = (\overline{\mathbb{T}}_+ \times M) \cup (\overline{\mathbb{T}}_- \times M)$. The Dixmier-Douady class trivializes in both components. Therefore, the groups that appear in the Mayer-Vietoris sequence are isomorphic to $G$-equivariant K-theory groups without a twist. Up to group isomorphisms, the Mayer-Vietoris sequence associated to the closed cover is (see \cite{TXL04} 3.12.) 
\begin{eqnarray*}
\xymatrix{ & K^0(\Gamma, \tau) \ar[r] &  K^0_G(\overline{\mathbb{T}}_+ \times M) \oplus K^0_G(\overline{\mathbb{T}}_- \times
M ) \ar[r]^{\hspace{1.2cm} a_0} &  K^0_G(\overline{\mathbb{T}}_{+-} \times M) \ar[d] \\ 
 & K^1_G(\overline{\mathbb{T}}_{+-} \times M) \ar[u] & \ar[l]_{\hspace{-1.2cm} a_1} K^1_G(\overline{\mathbb{T}}_+ \times
M) \oplus K^1_G(\overline{\mathbb{T}}_- \times M )  & \ar[l]  K^1(\Gamma,\tau) }.
\end{eqnarray*}
The maps $a_*$ are defined by $a_*(x,y) = (x-y, x - y \otimes \lambda_G)$ where $\lambda_G$ is the equivariant twisting line bundle. 

The solution of $K^*(\Gamma,\tau)$ from this sequence is parallel to \cite{HM12} and only the result is provided. One finds the solution by a diagram chase of an exact sequence of abelian groups. Therefore, the group structure is only recovered up to a possible extension problem. \3

\noindent \textbf{Proposition.}  If $* = 0,1$, then $K^{*}(\Gamma,\tau)$ is isomorphic to an extension of the group 
\begin{eqnarray*}
	\{x \in K^*_G(M): x = x \otimes \lambda_G \} \5 \text{by} \5 \frac{K^{*-1}_G(M)}{K^{*-1}_G(M) \otimes (1 - \lambda_G)}.
\end{eqnarray*}

\noindent \textbf{2.9.} Consider a $G$-equivariant projective Fock bundle $\textbf{PF}_G$ constructed in 2.6. We fix a $G$-equivariant finite rank complex vector bundle $\xi_G$ on $M$ and then define another projective bundle $\textbf{PF}_{\xi} = \textbf{P}(\textbf{F}_G \otimes \xi_G)$. Without a loss of generality, we can assume that $\lambda_G$ and $\xi_G$ are trivialized over the same cover $\mathfrak{U}$. We shall also apply the spinor module $\mathcal{H}_s$ and form a trivial bundle of projective spinors $\textbf{PS}$ as in 1.4. We have the projective bundle $\textbf{PS} \otimes \textbf{PF}_{\xi}$ which still has the Dixmier-Douady class represented by the cup product class $\tau$. 

The projective bundle $\textbf{PS} \otimes \textbf{PF}_{\xi}$ can be lifted locally, on each $U_i$ with $1 \leq i \leq 2n$, to a Hilbert bundle. Acting on the local Hilbert bundles we have the local Fredholm families $Q^i: U_i \rightarrow \textbf{Fred}_{\Psi}^{(1)}$ defined by
\begin{eqnarray*}
	Q^i_{(\phi,p)} = \sum_{n \in \mathbb{Z}} \psi_n \otimes e_{-n} + \frac{\phi}{2\pi} \psi_0 \otimes \textbf{1}, 
\end{eqnarray*}
The operators $e_n$ and $\psi_n$ are defined as in 1.4. Since the cover is $G$-invariant and the local families $Q^i$ are constant in the direction $M$, we get 
\begin{eqnarray*}	
g_{i,j+n}^{(1)}(a,x) Q^{j+n}_{(\phi,p)} (g_{i,j+n}^{(1)}(a^{-1}, ax))^{-1} = Q^i_{(\phi - 2 \pi,ap)} = Q^i_{a(\phi - 2 \pi,p)}.
\end{eqnarray*} 
for all $1 \leq i,j \leq n$ and $a \in G$. Consequely, $Q$ is a $G$-invariant section of a bundle of unbounded $\textbf{Fred}_{\Psi}^{(1)}$ operators associated to the projective bundle $\textbf{PS} \otimes \textbf{PF}_{\xi}$.\3

\noindent \textbf{Theorem.} The approximated sign of the section $Q$ defines a class in twisted K-theory:
\begin{eqnarray*}
F= \frac{Q}{\sqrt{1 + Q^2}} \in K^1_G(\mathbb{T} \times M, \tau).
\end{eqnarray*}

\noindent Proof. The $\Gamma$-invariance of $F$ follows immediately from the $\Gamma$-invariance of $Q$. The continuity of $F$ can be proved as in the case of smooth manifolds, \cite{HM12}. \5 $\square$ \3

Now the main problem to solve how the supercharge $Q$ associated with the vacuum bundle $\xi$ corresponds to the twisted K-theory group $K^1(\mathbb{T} \times M, \tau)$. An index-theoretic solution for this problem is given in 2.11. A sufficient background in equivariant cohomology is collected in 2.10, see \cite{BGV04}. \3

\noindent \textbf{2.10.} The groupoid cohomology groups $H^k(\Gamma, \mathbb{Z})$ are isomorphic to the $G$-equivariant integer cohomology groups because our groupoids are Morita equivalent to an action groupoid. We shall use the Cartan model of equivariant cohomology to study equivariant characters associated to the equivariant supercharges.   

Suppose that $N$ is a smooth manifold, and $G$ is a compact Lie group which acts smoothly on $N$. Let $\g$ denote the Lie algebra of $G$. If $X \in \g$ we denote by $X_N$ the vector field induced by the action $\exp(-tX)$ on $N$. Denote by $\iota_X : \Lambda^*(N) \rightarrow \Lambda^{* -1 }(N)$ the contraction by $X_N$. Define 
\begin{eqnarray*}
	\Lambda^q_G(N) =  \bigoplus_{2i + j = q} (S^i(\g^*) \otimes \Lambda^j(N))^G.
\end{eqnarray*}
This is the order $q$ subspace in the space of equivariant polynomial functions $\g \rightarrow \Lambda^{*}(N)$. The order zero subspace consists of $G$-equivariant functions. The equivariant coboundary $d_G: \Lambda^q_G(N) \rightarrow \Lambda^{q+1}_G(N)$ is defined by 
\begin{eqnarray*}
	(d_G \alpha)(X) = d(\alpha(X)) - \iota_X(\alpha(X))
\end{eqnarray*}
The square $d_G^2 \alpha(X) = - \{d, \iota_X \}\cdot \alpha(X) = - \mathcal{L}_X \cdot \alpha(X) $ is equal to zero since $\alpha(X)$ is invariant under the action of Lie derivative $\mathcal{L}_X$. We denote by $H^{*}_G(N, \mathbb{R})$ the cohomology of this complex. 

If $\nabla$ is an invariant connection of an equivariant line bundle $\lambda_G$, then the equivariant curvature of $\nabla$ is $F(X) = (\nabla - \iota_X)^2 + \mathcal{L}_X$. The first Chern class, $c^G_1$, sends an isomorphism class of a $G$-invariant line bundle to the cohomology class of its equivariant curvature form $(2 \pi i)^{-1} F(X)$. With this normalization the Chern class is an integer class in the equivariant cohomology group. The tensor product of equivariant line bundles defines a group structure in the set of isomorphism classes of $G$-equivariant line bundles. Then $c_1^G$ an isomorphism of groups. \3

\noindent \textbf{2.11.} We pull the projective bundle $\textbf{PS} \otimes \textbf{PF}_{\xi}$ and the supercharge $Q$ to the covering space $\mathbb{R} \times M$. The associated covering transformation is chosen so that it restricts to the identity over $M$. The Dixmier-Douady class becomes trivial on the covering space. The lifting of the $PU(\hil)$-valued transition functions to $U(\hil)$-valued operators on the covering space results a bundle of Hilbert space which is isomorphic to the trivial bundle $\mathbb{R} \times M \times PU(\hil)$. Let us denote this bundle by $\textbf{S} \otimes \textbf{F}_{\xi}$.

Now we can use standard methods in index theory to compute an equivariant characteristic polynomial on the covering space. Suppose that $\mathbb{A}$ is a $G$-invariant superconnection of the lifted supercharge. The $G$-invariance means that $\mathbb{A}$ commutes with the action. The supercurvature of $\mathbb{A}$ is given by the formula
\begin{eqnarray*}
\mathbb{F}(X) = (\mathbb{A} - \iota_X)^2 + \mathcal{L}_X
\end{eqnarray*}
and the Chern character $\mathbb{F}$ is the equivariant form
\begin{eqnarray*}
\textbf{ch-ind}_1(X) = \varphi (\text{sTr}(e^{- \mathbb{F}(X)})) \5 \text{for all } X \in \g,
\end{eqnarray*}
where sTr is the supertrace in the odd superconnection formalism and $\varphi$ is a normalization over $M$ which sends an equivariant $2n+1$ form $\Omega$ on $M$ to $(2 \pi i)^{-n} \Omega$. \3

\noindent \textbf{Proposition.}
\begin{quote}
\noindent (1) The character $\textbf{ch-ind}_1$ is a closed equivariant form.

\noindent (2) If $\mathbb{A}_t$ is a continuous one parameter family of superconnections with $t > 0$, then the transgression formula holds: 
\begin{eqnarray*}
	\frac{d}{dt} \textbf{ch-ind}_1(X) = - d_G \varphi(\text{sTr}(\frac{d \mathbb{F}_t(X)}{dt}e^{- \mathbb{F}_t(X)})). 
\end{eqnarray*}
\end{quote}

\noindent Proof. (1) The equivariant supercurvature is of the form $\mathbb{F}(X) = (d + \mathbb{B} - i_X)^2 + \mathcal{L}_X$ where $\mathbb{B}$ is a $G$-invariant endomorphism valued differential form which is odd with respect to the grading. Especially, $\text{sTr}[\mathbb{B}, \mathbb{F}(X)] = 0$ and one easily sees that $\text{sTr}[\mathbb{B}, e^{- \mathbb{F}(X)}] = 0$. Now
\begin{eqnarray*}
d_G  (\text{sTr} e^{-\mathbb{F}(X)}) = \text{sTr} [d_G, e^{-\mathbb{F}(X)}] = \text{sTr} [d_G + \mathbb{B}, e^{-\mathbb{F}(X)}] = 0.
\end{eqnarray*}
(2) The original proof, \cite{Qui85}, generalizes to the equivariant case. \3 
 
Let $\nabla_{\xi}$ and $\nabla_{\lambda}$ be $G$-invariant connections of the bundles $\xi_G$ and $\lambda_G$. A connection of $\textbf{S} \otimes \textbf{F}_{\xi}$ can be chosen by 
\begin{eqnarray*}
\nabla = N \nabla_{\lambda} + \nabla_{\xi}.
\end{eqnarray*}
Let $\hat{F}$ denote its curvature. Recall that $N$ is the unbounded operator which counts the charge of the Fock state. The strategy to compute the super Chern character is based on the homotopy invariance (2). We set a one parameter family of supercharges $t \mapsto \sqrt{t} Q = Q_t$ for $t > 0$ and define an equivariant superconnection by $\mathbb{A}_t = \chi Q_t + \nabla $. The equivariant supercurvature is
\begin{eqnarray*}
   \mathbb{F}_t(X) &=& (\chi Q_t + \nabla - \iota_X)^2 + \mathcal{L}_X \\
	&=& Q_t^2 + \{\chi Q_t, \nabla - \iota_X\} + \hat{F}(X) \\
	&=&Q_t^2 - \chi[\nabla, Q_t] + \hat{F}(X). 
\end{eqnarray*}
for all $X \in \g$. The formal symbol $\chi$ is applied in the odd superconnection formalism. It satisfies $\chi^2 = 1$, anticommutes $\iota_X$ and equivariant forms of odd rank and commutes with $Q$. The supertrace sTr applies the usual operator trace to the terms linear in $\chi$ and maps to zero everything else. The $t \rightarrow \infty$ limit of this character associated with this supercurvature can be computed in the equivariant case as in \cite{HM12}. The limit localizes as a distributional valued form  on the submanifold where $Q$ has its kernel. We let $P$ denote a family of projection operators onto the subbundle 
\begin{eqnarray*}
\eta_0 \otimes \bigoplus_{k \in \mathbb{Z}} S^k |0 \rangle \otimes \xi
\end{eqnarray*}
of the Hilbert bundle on the cover $\mathbb{R} \times M$.

Denote by $y \in \mathbb{R}$ the lifted coordinate. \3

\noindent \textbf{Proposition.} The $t \rightarrow \infty$ limit of the form $\textbf{ch-ind}_1$ is the distribution valued form 
\begin{eqnarray*}
	\lim_{t \rightarrow \infty} \textbf{ch-ind}_1 = \sqrt{\pi}  \delta(Pe_0P + y) dy \wedge \ch_G(\hat{F})(X), 
\end{eqnarray*}
where $\delta$ denotes the Dirac delta distribution and $\ch(\xi)(X)$ is the equivariant Chern character of the connection $\hat{F}$: 
\begin{eqnarray*}
\ch_G(\hat{F})(X) = \tr_{\xi} e^{- \frac{\hat{F}(X)}{2\pi i}}. 
\end{eqnarray*}

We would like to fix a locally defined section $\psi: \mathbb{T} \times M \rightarrow \mathbb{R} \times M$ and pull the form $\textbf{ch-ind}_1(X)$ to the base $\mathbb{T} \times M$. There is the difficulty that under translations by $2 \pi$ in the direction of $\mathbb{R}$, the supercurvature transforms as 
\begin{eqnarray*}
\mathbb{F}_t(X)_{(y,p)} &\mapsto& S (\mathbb{F}_{t}(X))_{(y + 2 \pi, p)} S^{-1} \\
&=& S(Q_t - \chi[\nabla, Q_t] + \hat{F}(X))_{(y + 2 \pi, p)}S^{-1} \\ 
&=& S(Q_t - \chi d y + \hat{F}(X))_{(y + 2 \pi, p)}S^{-1} \\
&=& (Q_t - \chi dy + \hat{F}(X) - F_{\lambda}(X))_{(y, p)} \\
&=& \mathbb{F}_t(X)_{(y,p)} - F_{\lambda}(X)_{(y, p)}
\end{eqnarray*}
for all $(y, p) \in \mathbb{R} \times M$. Therefore, the superconnection is not invariant under translation by $2 \pi$ and it would not glue to a global form under the pullback. 

Given a section $\psi$ of the covering map $\mathbb{R} \times M \rightarrow \mathbb{T} \times M$, we use the strategy in \cite{HM12} and pull the form to a quotient of the rational equivariant cohomology group
\begin{eqnarray*}
\underline{\textbf{ch-ind}}^{\tau}_1 = \psi^*(\textbf{ch-ind}_1) \in \frac{\sqrt{\pi} H^{\mathrm{odd}}_G(\mathbb{T} \times M, \mathbb{Q})}{\sqrt{\pi} \frac{d\phi}{2 \pi} \wedge (1 - \ch_G(F_{\lambda}))\wedge \ch_G(K_G(M))}.
\end{eqnarray*}
Here we consider $H^{\mathrm{odd}}_G(\mathbb{T} \times M, \mathbb{Q})$ to be an abelian group and do not equip it with a $\mathbb{Q}$-module structure. Likewise, $\sqrt{\pi} \frac{d\phi}{2 \pi} \wedge (1 - \ch_G(F_{\lambda}))\wedge \ch(K_G(M))$ is the normal subgroup of cohomology classes of these elements under addition. In conclusion. \3

\noindent \textbf{Theorem.} Let $\psi$ be an arbitrary section of the covering map. The equivalence class of the equivariant character $\underline{\textbf{ch-ind}}^{\tau}_1$ is independent on the choice of the section $\psi$ and is represented by 
\begin{eqnarray*}
\sqrt{\pi}  \frac{d \phi}{2 \pi} \wedge \ch_G(F_{\xi}) \in \frac{\sqrt{\pi} H^{\mathrm{odd}}_G(\mathbb{T} \times M, \mathbb{Q})}{\sqrt{\pi}\frac{d\phi}{2 \pi} \wedge (1 - \ch_G(F_{\lambda}))\wedge \ch_G(K_G(M))}.
\end{eqnarray*}

\noindent Proof. The distributional form can be represented in cohomology by the angular form with volume equal to $1$. Thus, the character can be represented by 
\begin{eqnarray*}
\underline{\textbf{ch-ind}}^{\tau}_1 = \sqrt{\pi} \psi^*( \frac{dy}{2 \pi} \wedge \ch_G(F_{\xi})). 
\end{eqnarray*}
Under the translation by $2 \pi$ in $\mathbb{R}$, the supercurvature will be shifted by $- F_{\lambda}$. Since $F_{\lambda}$ is $G$-invariant and thus commutes with the supercurvature, the character form gets multiplied by $\ch(F_{\lambda})$. However, this operation does not have any effect in the quotient where the character is evaluated. Therefore, the choice of $\psi$ does not have effect on the equivalence class. \5 $\square$ \3

This index theorem can be used to extract information from the twisted K-theory elements represented by approximated signs of the supercharges. Namely, if $Q_1$ and $Q_2$ is a pair of supercharges with vacuum bundles $\xi_1$ and $\xi_2$, and if the characteristic classes of $Q_1$ and $Q_2$ are not equal, then there can be no continuous homotopy between $Q_1$ and $Q_2$. Since the topology in the space of unbounded Fredholm operators is induced from the topology of bounded Fredholm operators, there can be no continuous homotopy at the level of approximated signs either. Thus, the twisted K-theory classes are different as well.   
 
\section{Examples} 

\noindent \textbf{3.1.} In the first  part of this series, \cite{HM12}, we gave explicit construction for the twisted $K^1$-groups on the product manifolds $\mathbb{T} \times M$ when $M = \mathbb{T}^2$ and $S^2$. In both cases, the usual $K$-theory groups can be computed, for example, using the Kunneth formula which leads to 
\begin{eqnarray*}
	K^1(\mathbb{T} \times S^2) \simeq \mathbb{Z}^{\oplus 2}, \5 K^1(\mathbb{T} \times \mathbb{T}^2) \simeq \mathbb{Z}^{\oplus 4}.
\end{eqnarray*}
The twisting line bundles are classified by integers. For $\lambda$ corresponding to $0 \neq k \in \mathbb{Z}$ the twisted $K^1$-groups are  
\begin{eqnarray*}
	K^1(\mathbb{T} \times S^2) \simeq \mathbb{Z} \oplus \mathbb{Z}_k, \5 K^1(\mathbb{T}^3) = \mathbb{Z}^{\oplus 3} \oplus \mathbb{Z}_k. 
\end{eqnarray*}
The torsion sugroups appear in the twisted theory.  In the case of real projective spaces one has torsion subgroups in the ordinary K-theory but they disappear after twisting. This is explained in 3.2. in the supercharge formalism. \3

\noindent \textbf{3.2.} We study the real projective space $\mathbb{R}P^{2n}$ for $n \in \mathbb{N}$. No group action on $\mathbb{R}P^{2n}$ is assumed here. The second cohomology group of $\mathbb{R}P^{2n}$ is the torsion group $\mathbb{Z}_2$. Let us denote by $\lambda$ the nontrivial line bundle which is subject to the relation $\lambda \otimes \lambda = 1$. The cohomology groups of $\mathbb{T} \times \mathbb{R}P^{2n}$ can by computed with Kunneth's formula. The third cohomology is  
\begin{eqnarray*}
	H^3(\mathbb{T} \times \mathbb{R}P^{2n}) = \mathbb{Z}_2 \5 \text{for all $n \in \mathbb{N}$.}
\end{eqnarray*}
There is exactly one nontrivial cup product twisting class $\tau = \alpha \smile \beta$ where $\beta$ is the class of $\lambda$ in $H^2(\mathbb{R}P^2)$ and $\alpha$ is the generator of $H^1(\mathbb{T})$. According to the proposition of 1.3, the projective Fock bundle model applied to the twisting line bundle $\lambda$ realizes this class geometrically.  

K-theory groups on $\mathbb{R}P^{2n}$ are solved in \cite{Ati67}
\begin{eqnarray*}
	K^0(\mathbb{R}P^{2n}) \simeq \mathbb{Z} \oplus \mathbb{Z}_{2^n}, \5 K^1(\mathbb{R}P^{2n}) = 0. 
\end{eqnarray*}
As an abelian group, $K^0$ is generated by $1$ and $x_n$ and its ring structure is given by  
\begin{eqnarray}\label{rpn}
	x_n = 1 - [\lambda], \5 x_n^2 = 2 x_n, \5 2^n x_n = 0.
\end{eqnarray}
By Kunneth's formula we have
\begin{eqnarray*}
	K^1(\mathbb{T} \times \mathbb{R}P^{2n}) \simeq \mathbb{Z} \oplus \mathbb{Z}_{2^n}.
\end{eqnarray*}
Next we use the proposition of 2.8 to solve the twisted $K^1$-group. The normal subgroup 
\begin{eqnarray*}
K^0(\mathbb{R}P^{2n})\otimes (1 - \lambda) = K^0(\mathbb{R}P^{2n})\otimes x_n
\end{eqnarray*}
is the torsion subgroup isomorphic to $\mathbb{Z}_{2^n}$ in $K^0(\mathbb{R}P^{2n})$ by \eqref{rpn}. Therefore, the $\tau$-twisted $K^1$-group is given by 
\begin{eqnarray*}
	K^1(\mathbb{T} \times \mathbb{R}P^{2n}, \tau) = \mathbb{Z}.
\end{eqnarray*}
The twisting by $\tau$ kills the torsion subgroup.

Consider a pullback of the projective bundle of Fock spaces $\textbf{PF}_{\xi}$ to $\mathbb{R} \times \mathbb{R}P^{2n}$ as in 2.11. Then we can lift the transition functions and define a Hilbert bundle on the covering space. Assume that $\xi$ is a line bundle on $\mathbb{R}P^{2n}$. Thus, $\xi = \lambda$ or $\xi = 1$, the trivial complex line bundle. Then the topological type of the vacuum at the submanifold with $0 = y \in \mathbb{R}$ is ether $\lambda$ or $1$. The type depends on the chosen lift $\mathbb{R} \rightarrow \mathbb{T}$ but is not important.  Translations in the direction of positive real numbers by $2 \pi$ then change this type to $1$ or $\lambda$ (respectively). This continues through the real line: the topological type of the vacuum is $\mathbb{Z}_2$ graded and depends on the choice of integer part of the $\mathbb{R}$ coordinate. This is the reason for the absence of the torsion subgroup: twisted $K^1$-theory does not distinguish $\lambda$ from $1$ since vacuums of these types are identified.  

In the nontwisted case (i.e. set $\beta = 0$), we have the supercharge construction on $\mathbb{T} \times \mathbb{R}P^{2n}$. Without the twist, the topological type of the Fock vacuum does not depend on the direction $\mathbb{T}$ and hence, the topological type of the kernel bundle is fixed by the choice of $\xi$. 

Let $\tau$ be the nontrivial twisting class. Now twisted K-theory is insensitive to the topology associated with $\lambda$ in the local kernel bundles. Thus, only the information associated with the rank of the kernel bundle is relevant. Suppose that the vacuum bundle $\xi$ is an arbitrary complex vector bundle. We can study the characteristic form of the supercharge $Q$ applying the methods of 2.11 with a trivial group action $G$. Detailed calculations can also be found from Section 6 in \cite{HM12}. The curvature form of $\lambda$ vanishes in the de Rham cohomology. Now the form $\underline{\textbf{ch-ind}}^{\tau}_1$ gets values in the usual odd de Rham cohomology of $\mathbb{T} \times \mathbb{R}P^{2n}$. In fact, according to the analysis in 2.11 the characteristi form is represented by  
\begin{eqnarray*}
\underline{\textbf{ch-ind}}^{\tau}_1 = \text{rank}(\xi) \sqrt{\pi} \frac{d\phi}{2 \pi} \in \sqrt{\pi} H^{\text{odd}}(\mathbb{T} \times \mathbb{R}P^{2n}, \mathbb{Q}). 
\end{eqnarray*}
The rank of the vacuum bundle determines the twisted $K^1$-class. \3
 
\noindent \textbf{3.3.} Let $\mathbb{T}$ act on $S^2$ by rotations around the axis of the sphere. We denote by $S^2_+$ and $S^2_-$ the north and south hemispheres of $S^2$ such that the intersection $S^2_{+-}$ is the equator $\mathbb{T}$. Then $S^2 = S^2_+ \cup S^2_-$ and the cover is $\mathbb{T}$-invariant. The south and north poles are the fixed points of this action. The equivariant line bundles are classified by 
\begin{eqnarray*}
	H^2_{\mathbb{T}}(S^2, \mathbb{Z}) \simeq \mathbb{Z}[u_+] \oplus \mathbb{Z}[u_-].
\end{eqnarray*}
This result is well known and can be solved with the Mayer-Vietoris sequence. $u_{\pm}$ are generators of the cohomology groups $H^2_{\mathbb{T}}(S^2_{\pm}, \mathbb{Z})$: $S^2_{\pm}$ are retracted to the fixed points resulting cohomology groups isomorphic to the $\mathbb{T}$-equivariant cohomology of a point, i.e. $H^2(\mathbb{C}P^{\infty}, \mathbb{Z})$, and under this identification, $u_{\pm}$ can be viewed as generators of $H^2(S^2, \mathbb{Z})$. We also get the isomorphism 
\begin{eqnarray*}
	H^3_{\mathbb{T}}(\mathbb{T} \times S^2,\mathbb{Z}) \simeq H^2_{\mathbb{T}}(S^2, \mathbb{Z})
\end{eqnarray*}
from an Mayer-Vietoris sequence. 

The ordinary equivariant K-theory can be computed using the Mayer-Vietoris sequence
\begin{eqnarray*}
\xymatrix{ & K_{\mathbb{T}}^0(S^2) \ar[r] &  K_{\mathbb{T}}^0(S^2_+) \oplus K_{\mathbb{T}}^0(S^2_-) \ar[r]^{\hspace{1.2cm} a_0} &  K_{\mathbb{T}}^0(S^2_{+-}) \ar[d]  \\ 
 & K_{\mathbb{T}}^1(S^2_{+-}) \ar[u] & \ar[l]_{\hspace{-1.2cm} a_1} K_{\mathbb{T}}^1(S^2_+) \oplus K_{\mathbb{T}}^1(S^2_-) & \ar[l] K_{\mathbb{T}}^1(S^2) }.
\end{eqnarray*}
The groups $K_{\mathbb{T}}^0(S^2_+)$ and $K_{\mathbb{T}}^0(S^2_-)$ are isomorphic to  $K^0_{\mathbb{T}}(*)$ where $*$ is a pole and $K^0_{\mathbb{T}}(*)$ is isomorphic to the representation ring
\begin{eqnarray*}
   R(\mathbb{T}) = \mathbb{Z}[a,a^{-1}]. 
\end{eqnarray*}
The circle group action on the equator $S^2_{+-}$ is free. Therefore
\begin{eqnarray*}
	&& K_{\mathbb{T}}^0(S^2_{+-}) \simeq K^0(\mathbb{T}/\mathbb{T}) \simeq \mathbb{Z} \5 \text{and} \\
	&& K_{\mathbb{T}}^1(S^2_{+-}) \simeq K^1(\mathbb{T}/\mathbb{T}) \simeq 0
\end{eqnarray*}
Next we observe that $K_{\mathbb{T}}^1(S^2_{+}) \simeq K_{\mathbb{T}}^1(*)$. Thus, $K_{\mathbb{T}}^1(S^2_{+}) = 0$ and similarly $K_{\mathbb{T}}^1(S^2_{-}) = 0$. The homomorphisms $a_i$, $i = 0,1$ send $(x,y)$ to $x-y$. The group elements in $K^0_{\mathbb{T}}(S^2_{+-})$ are determined by an integer, the virtual rank. Thus, 
\begin{eqnarray*}
	\text{Ker}(a_0) \simeq \frac{R(\mathbb{T}) \oplus R(\mathbb{T})}{ \mathbb{Z}}, \5 \text{Im}(a_0) \simeq \mathbb{Z} \5 \text{and} \5 \text{Ker}(a_1) = 0.
\end{eqnarray*}
Since there is no torsion in these groups the K-theory groups are as follows 
\begin{eqnarray*}
	&& K^0_{\mathbb{T}}(S^2) \simeq \frac{K^1_{\mathbb{T}}(S^2_{+-})}{\mathrm{Im}(a_1)} \oplus \mathrm{Ker}(a_0) \simeq  \frac{R(\mathbb{T}) \oplus R(\mathbb{T})}{\mathbb{Z}} \\
	&& K^1_{\mathbb{T}}(S^2) \simeq \frac{K^0_{\mathbb{T}}(S^2_{+-})}{\mathrm{Im}(a_0)} \oplus \mathrm{Ker}(a_1) \simeq 0.
\end{eqnarray*}
The elements in $K^0_{\mathbb{T}}(S^2)$ are the pairs $(x,y) \in R(\mathbb{T})^{\oplus 2}$ constrained by the common virtual dimension. Topologically any vector bundle over $S^2$ is a direct summand of line bundles. Therefore, $K^0_{\mathbb{T}}(S^2)$ is the abelian group generated by the $\mathbb{T}$-equivariant complex line bundles and each generator carries a representations $(a^{n_+}, a^{n_-}) \in R(\mathbb{T})^{\oplus 2}$ on its fibres. Then $\mathbb{T}$ acts on the fibers over $S^2_{\pm}$ through the character $\exp(i \varphi) \mapsto \exp(i n_{\pm}\varphi)$. 

The twisted $K^1$-group can be solved with the proposition of 2.8. If $\lambda_{\mathbb{T}}$ is an equivariant line bundle associated with the reprentation $(a^{l_+},a^{l_-})$ and the $K^1$-group is subject to the twisting of type $\alpha \smile \beta_{\mathbb{T}}$ and $\beta_{\mathbb{T}}$ is the 2-cohomology class which identifies with the equivalence class of $\lambda_{\mathbb{T}}$, then we have 
\begin{eqnarray*}
	K^1_{\mathbb{T}}(\mathbb{T} \times S^2, \tau) \simeq \frac{K^0_{\mathbb{T}}(S^2)}{K^0_{\mathbb{T}}(S^2) \otimes (1 - \lambda_{\mathbb{T}})}.
\end{eqnarray*}
Therefore, the twisted $K^1$-group is isomorphic to the abelian group defined by 
\begin{eqnarray*}
	K^1_{\mathbb{T}}(\mathbb{T} \times S^2, \tau) \simeq K^0_{\mathbb{T}}(S^2)/ \sim \5 \text{where} \5 (a^{n_+}, a^{n_-}) \sim (a^{n_+ + l_+}, a^{n_- + l_-}).
\end{eqnarray*} 

We can find explicit groupoid cocycles for the $\mathbb{T}$-equivariant gerbes on $\mathbb{T} \times S^2$. Now we replace the closed sets $S^2_{\pm}$ with $\mathbb{T}$-invariant open sets which cover the north and the south hemispheres and which intersect in a tube around the equator. These sets are still denoted by $S^2_{\pm}$. Then there is a Lie groupoid
\begin{eqnarray*}
\Xi = \coprod_{ij} \mathbb{T} \times S^2_{ij} \rightrightarrows \coprod_i S^2_{i}
\end{eqnarray*}
where $i,j \in \{+,-\}$. This is Morita equivalent to the action groupoid $\mathbb{T} \ltimes S^2 \rightrightarrows S^2$. The equivariant complex line bundle corresponding to $(a^{n_+}, a^{n_-})$ in $H^2_{\mathbb{T}}(S^2)$ is associated to the groupoid cocycle $h \in H^1(\Xi, \underline{\mathbb{T}})$. Its components are
\begin{eqnarray*}
	&& h_{++}(\exp(i \varphi), p) = \exp(i n_+ \varphi),\\
	&& h_{--}(\exp(i \varphi), p) = \exp(i n_- \varphi),\\
	&& h_{+-}(\exp(i \varphi), (\alpha,r)) = \exp(i n_+ \varphi + i(n_+ - n_-) \alpha) \5 \text{and} \\
	&& h_{-+}(\exp(i \varphi), (\alpha,r)) = \exp(i n_- \varphi - i(n_+ - n_-) \alpha)
\end{eqnarray*}
where we have denoted by $\alpha$ the angle coordinate and $r$ the height coordinate of the tube around the equator $S^2_{+-}$. The components $h_{++}$ and $h_{--}$ carry the information of the $\mathbb{T}$-action on the fibres. The topological type of this bundle is determined by $h_{+-}(1,(\alpha,r))$ which is a function on the equator with winding number equal to $(n_+ - n_-)$. 

Denote by $\{U_i\}$ the cover of $\mathbb{T} \times S^2$ which has 4-components: $\{\mathbb{T}_{\pm} \times S^2_{\pm} \}$. Following 2.4 we construct the Lie groupoid
\begin{eqnarray*}
\Gamma = \mathbb{T} \times \coprod_{ij} U_{ij} \rightrightarrows \coprod_i U_i.
\end{eqnarray*}
The cocycle $g \in H^1(\Gamma, \underline{PU}(\hil))$ of the projective Fock bundle is straightforward to write down using the formulas in 2.6. 
 
Next we tensor the Fock vacuum with $\xi_{\mathbb{T}}$ which is a line bundle corresponding to the class $(a^{j_+}, a^{j_-})$. Let $F_{\lambda}$ and $F_{\xi}$ denote $G$-invariant curvatures of $\lambda_{\mathbb{T}}$ and $\xi_{\mathbb{T}}$. We apply the theorem of 2.11 which gives a representative for the character map of the supercharge associated with this setup. The 3-form part reads 
\begin{eqnarray*}
	- \sqrt{\pi} \frac{d\phi}{2 \pi} \wedge \frac{F_{\xi}}{2 \pi i} \5 &\mathrm{mod}& \5 \sqrt{\pi} \frac{d\phi}{2 \pi} \wedge \frac{F_{\lambda}}{2 \pi i} \5 \Leftrightarrow \\
	-\sqrt{\pi} \frac{d\phi}{2 \pi} \wedge ( j^+ F^+ + j^- F^-) \5 &\text{mod}& \5  \sqrt{\pi} \frac{d\phi}{2 \pi} \wedge ( l^+ F^+ + l^- F^-).
\end{eqnarray*}
The equivariant forms $F^{\pm}$ are generators of $H^2_{\mathbb{T}}(S^2, \mathbb{Z})$. Therefore, if $\xi_1$ and $\xi_2$ are equivariant line bundles associated to the classes $(a^{j^+_1}, a^{j^-_1})$ and  $(a^{j^+_2}, a^{j^-_2})$, and $\lambda_{\mathbb{T}}$ is as above, then the supercharges associated with these models determine distinct elements in $K^1_{\mathbb{T}}(\mathbb{T} \times S^2, \tau)$ if and only if
\begin{eqnarray*}
	j^+_1 - j^+_2 \neq 0 \5 \text{mod}\5 l^+ \5 \text{or} \5 j^-_1 - j^-_2 \neq 0 \5 \text{mod}\5 l^-.
\end{eqnarray*}

If we apply similar analysis as above for higher rank $\mathbb{T}$-equivariant vector budles $\xi_1$ and $\xi_2$ the characteristic map is no longer injective as a map from supercharges to equivariant cohomology theory. However, the superconnection analysis still provides nontrivial information. For example, if the ranks of $\xi_1$ and $\xi_2$ do not match, then the associated $K^1$-classes are necessarily different since the one form components of the characters are not equal in cohomology. 

\section*{Appendix}

\noindent \textbf{A.1.} This appendix is motivated by the geometric and group theoretic structures in Hamiltonian quantization. In quantum field theory on a unit circle $\mathbb{T}$, one defines a space of potentials $\mathcal{A}$ to be the affine space of $\g$-valued one forms on $\mathbb{T}$ and $\g$ is a Lie algebra of a simple, compact and simply connected Lie group $G$. Then there is a family of Dirac operators acting on  $L^2(\mathbb{T})$ which are coupled to the potentials. The gauge symmetry group is chosen to be the group of based gauge transformations $\Omega G$ which is the subgroup in the group of smooth loops $LG$ based at the unit $\gamma(0) = 1 = \gamma(2\pi)$. The true parameter space of the theory is the quotient $\mathcal{A}/\Omega G$ where $\Omega G$ acts on $\mathcal{A}$ by
\begin{eqnarray*}
	(g,A) \mapsto g A g^{-1} + g d(g^{-1}).
\end{eqnarray*}
The action makes $\mathcal{A} \rightarrow \mathcal{A}/\Omega G$ a principal bundle with a projection map sending $A \in \mathcal{A}$ to a holonomy around the circle.

On $\mathcal{A}/\Omega G$ which is homeomorphic to $G$ we define a Fock bundle. We set vacuum levels $x_i$ and a finite cover $U_i$ of $\mathcal{A}/\Omega G$ so that the the numbers $x_i$ do not belong to the spectrum of the Dirac family over $U_i$. Then, over each $U_i$, one can define a polarization in the local Hilbert bundles to Dirac eigenspaces with eigenvalues larger than $x_i$ and its orthogonal complement. The global spectral flow of the Dirac family is an obstruction to set this polarizations globally over $\mathcal{A}/\Omega G$. Next one fixes local Fock bundles associated with the chosen polarization on each $U_i$. On the overlaps $U_{ij}$ two different vacuums associated to different polarizations are related by tensoring with  a determinant line bundle. The Hamiltonian system is therefore naturally defines as a projective bundle of Fock spaces. Then there is a prolongation problem: can we fix the vacuum levels locally (local determinant line bundles) to create a Fock bundle over $\mathcal{A}/\Omega G$. In \cite{CMM97}, \cite{Lot02} it was found that the obstruction to do this is the three cohomology part of the odd character of the Dirac family over $\mathcal{A}/ \Omega G$. The gerbe analysis in the first part of this series \cite{HM12} was based on this result. 

The gauge group $\Omega G$ acts on the fibers of the Fock bundle as a subgroup of $LG$. The loop group acts under an irredcible positive energy representation which is projective. Alternatively, these representations are irreducible representations of the central extension $\hat{LG}$. The extension cocycles arise in the lifted transition functions of the Fock bundle over $\mathcal{A} / \Omega G$ and ruin the {\v C}ech cocycle condition for the transition functions. The relationship to the prolongation problem can be described explicitly.  Let $c \in H^2(\mathfrak{lg}, \mathbb{C})$ denote the Lie algebra cohomology class of the extension. We can identify the Lie algebra elements with the tangent vector fields near the identity and realize $c$ as a differential two form near the identity of $LG$ and extend it by left translation through $LG$. The normalized cocycle  $c / 2 \pi i $ defines an integer cohomology class. We choose local sections $s_i: U_i \rightarrow p^{-1}(U_i)$ of $\mathcal{A} \stackrel{p}{\rightarrow} \mathcal{A}/\Omega G$. Let $g_{ij}$ denote the $\Omega G$ valued transition functions of the principal bundle $\mathcal{A}$. Then $s_j = s_i g_{ij}$. We pull back the two forms $c$ on $LG$ to define local 2-forms $g_{ij}^*(c/2 \pi i)$ on $U_{ij}$. These define a representative for the Dixmier-Douady class of the gerbe in {\v C}ech-de Rham double complex. The {\v C}ech-de Rham complex can be analyzed in as usual to get a pure de Rham representative. We find the forms $\omega_i$ and $\omega_j$ over $U_i$ and $U_j$ so that $g^*_{ij}(c/2 \pi i) = \omega_j - \omega_i$. Let $\Omega$ be a 3-form over $G$ which is locally defined by $d \omega_i$ on $U_i$. Then $\Omega$ is a global 3-form which represents the Dixmier-Douady class of the gerbe. Explicitly, 
\begin{eqnarray*}
	\Omega = k\frac{1}{24 \pi^2}(g^{-1} dg)^3
\end{eqnarray*}
is $k$ times the generator of $H^3(G, \mathbb{Z})$ if the projective representation of $LG$ has level equal to $k \in \mathbb{Z}$. For details, see \cite{Mic06}, \cite{MPel07}. \3

\noindent \textbf{A.2.} We can view the torus $\mathbb{T}^3$ as a parameter space for a quantum field theory in $\mathbb{T}_{\theta}$: we replace $\mathcal{A}$ by $\mathbb{R}^3$ and $\Omega G$ by $\mathbb{Z}^3$. The action is by translations and therefore $\mathbb{R}^3 / \mathbb{Z}^3 = \mathbb{T}^3$. We write $\mathbb{T}^3 = \mathbb{T}_{\phi} \times \mathbb{T}^2$. In \cite{HM12} we constructed a Hilbert bundle over the 4-torus $\mathbb{T}_{\theta} \times \mathbb{T}^3$ with fibers $L^2(\mathbb{T}_{\theta})$. There is a Dirac family over $\mathbb{T}_{\theta} \times \mathbb{T}^3$ coupled to the connection of the bundle $\lambda_1 \boxtimes \lambda_k$ where $\lambda_i$ is the line bundle associated to the Chern class $i \in \mathbb{Z}$ on a 2-torus. The 3 cohomology component of the character of this Dirac family is the cup product $\alpha \smile \beta$ such that $\alpha$ is the generator of $H^1_{\mathrm{dR}}(\mathbb{T}_{\phi}, \mathbb{Z})$ and $\beta$ is  $(2 \pi i )^{-1}$ times the curvature form of $\lambda_k$. The associated decomposable {\v C}ech cocycle defines a transition data of the gerbe exactly as in 1.3.

Here we want to connect the topological prolongation problem to a group extension problem: we can define a trivial bundle of Fock spaces on the affine space $\mathbb{R}^3$ and then implement an extended $\mathbb{Z}^3$-symmetry on the fibres by an extension cocycle which defines the cup product class in cohomology. In the following we apply the theory of groupoid extensions to make the connection between this extension problem and the cohomology class of the gerbe explicit.\3

\noindent \textbf{A.3.} We can consider the torus $\mathbb{T}^3$ as a transformation groupoid $\Gamma := \mathbb{Z}^3 \times \mathbb{R}^3  \rightrightarrows \mathbb{R}^3$ with the usual target and source maps 
\begin{eqnarray*}
t,s: \Gamma \rightrightarrows \mathbb{R}^3, \5 s(a,x) = x, \5 t(a,x) = x+a. 
\end{eqnarray*}
The composition $m: \Gamma^{(2)} \rightarrow \Gamma$ is defined in its domain by $(b,x+a)(a,x) = (a+b,x)$. The action is free and therefore the cohomology of this groupoid is just the usual de Rham cohomology of the 3-torus. 

Next we define a topologically trivial $\mathbb{T}$-extension $\hat{\Gamma}:= \mathbb{Z}^3 \times \mathbb{R}^3 \times \mathbb{T} \rightrightarrows \mathbb{R}^3$ with a projection $p: \hat{\Gamma} \rightarrow \Gamma$ defined by $p(a,x,\mu) = (a,x)$. The target and source maps of the groupoid $\hat{\Gamma}$ are 
\begin{eqnarray*}
\hat{s} = s \circ p, \5 \hat{t} = t \circ p.
\end{eqnarray*}
The nontriviality of the extension is associated to the cocycle $\omega \in H^2(\Gamma, \underline{\mathbb{T}})$ which is represented by a map $\omega: \Gamma^{(2)} \rightarrow \mathbb{T}$ and the cocycle condition reads
\begin{eqnarray*}
 (\partial^* (\omega))((a,x + b + c), (b,x + c), (c,x)) = 0
\end{eqnarray*}
which is equivalent to 
\begin{eqnarray*}
	& & \omega((b,x + c),(c,x)) \omega((a,x + b + c),(b+c,x)) = \\
	& & \omega((a+b,x+c),(c,x)) \omega((a,x+b+c),(b,x+c)).
\end{eqnarray*}
We set $x = (x_0, x_1, x_2)$ and $a = (a_0, a_1, a_2)$. One checks easily that $\omega: \Gamma^{(2)} \rightarrow \mathbb{T}$ defined by 
\begin{eqnarray}\label{t3}
\omega((a,x+b),(b,x)) = e^{2\pi i k a_0 b_2 x_1}.
\end{eqnarray}
is a cocycle. Then we define a composition in the extension groupoid $\hat{\Gamma}$ by 
\begin{eqnarray*}
(a,x + b,\mu')(b,x,\mu) = (a+b,x,  \mu \mu' \omega((a,x+b),(b,x))).
\end{eqnarray*}
This extensions corresponds to the class $k \in \mathbb{Z} \simeq H^3(\Gamma, \mathbb{Z})$. This will be justified below.
 
One can define generators $\widetilde{\alpha}_i$,  $i = 0,1,2$ for $H^1(\Gamma, \mathbb{Z}) \simeq \mathbb{Z}^3$ by
\begin{eqnarray*}
	\widetilde{\alpha}_i(a,x) = a_i \5 \hbox{for all}\5 (a,x) \in \mathbb{Z}^3 \times \mathbb{R}^3. 
\end{eqnarray*}
The cohomology classes of a 3-torus can be constructed by applying the cup product to the generators. Then the 3-cohomology group is generated by the class of
\begin{eqnarray*}
	(\widetilde{\alpha}_2 \smile \widetilde{\alpha}_1 \smile \widetilde{\alpha}_0) ((a, x + b + c), (b,x + c), (c,x)) = a_2 b_1 c_0. 
\end{eqnarray*}
We choose $\widetilde{\beta}$ to be the cohomology class represented by $k \widetilde{\alpha}_2 \smile \widetilde{\alpha}_1$. Then we write $\widetilde{\alpha} = \widetilde{\alpha}_0$ and
\begin{eqnarray*}
	(\widetilde{\alpha} \smile \widetilde{\beta})((a, x + b + c), (b,x + c), (c,x)) =  k a_0  b_2 c_1.
\end{eqnarray*}
By antisymmetry, this cocycle is equivalent to $k$ times the generator. Next we define a cocycle $f' \in H^2(\Gamma, \underline{i\mathbb{R}})$ by 
\begin{eqnarray*}
	f'((a,x+b),(b,x)) =  2 \pi i k a_0 b_2 x_1.
\end{eqnarray*}
It satisfies $(2 \pi i)^{-1} \partial^*(f') = \widetilde{\alpha} \smile \widetilde{\beta}$. It then follows that under the usual isomorphism, $\widetilde{\alpha} \smile \widetilde{\beta}$ maps to the cocycle $\omega \in H^2(\Gamma, \underline{\mathbb{T}})$ defined in \eqref{t3}.

In the action groupoid formalism there is a globally trivial $\mathbb{Z}^3$-equivariant projective Hilbert bundle over $\mathbb{R}^3$ with a $H^1(\Gamma, \underline{PU}(\mathcal{H}))$ class uniquely determined (in cohomology) by $\omega$ in \eqref{t3}. The equivariant bundle can be chosen to be topologically a product $\mathbb{R}^3 \times \mathcal{PF}$ where $\mathcal{PF}$ is a projective Fock space. Given a cocycle  $g \in H^1(\Gamma, \underline{PU}(\mathcal{H}))$, the fibres transform under the action of $a \in  \mathbb{Z}^3$ by
\begin{eqnarray*}
	(x, \Psi) \mapsto (x+a, g(a,x) \Psi).
\end{eqnarray*}
The cocycle condition $\partial^*(g) = 0$ is equivalent to the associativity of the action. If we lift the action to the group $U(\mathcal{H})$ then $\omega \in H^2(\Gamma, \underline{\mathbb{T}})$ is risen from the relation
\begin{eqnarray*}
	\partial^*(\hat{g}) = \omega \5 \Leftrightarrow \5 \hat{g}(a,x+b)\hat{g}(b,x)\hat{g}(a+b,x)^{-1} = \omega((a,x+b),(b,x)) 
\end{eqnarray*}
where $\hat{g}$ are the lifted functions $\Gamma^{(1)} \rightarrow U(\mathcal{H})$. 

We define the cocycle $g\in H^1(\Gamma, \underline{PU}(\hil))$ by ($S$ and $N$ are defined as in 1.3) 
\begin{eqnarray}\label{puh}
	g(a,x) = (e^{2 \pi i k a_2  x_1})^N S^{a_0}.
\end{eqnarray}
The cocycle condition is straightforward to check. When lifted to $U(\hil)$, one finds a representative for the Dixmier-Douady class
\begin{eqnarray*}
	\omega((a,x+b),(b,x))  = e^{2 \pi i k a_0 b_2 x_1}. 
\end{eqnarray*}
This is the class \eqref{t3}. 

If we go to the quotient space $\mathbb{T}^3$, then the action \eqref{puh} on the fibres defines a gerbe of topological type $\alpha \smile \beta$ where $\beta$ is a {\v C}ech cohomology class associated to the line bundle with Chern class equal to $k$ and $\alpha$ is the {\v C}ech generator of $\mathbb{T}$.

\end{document}